\documentclass[11pt]{article}
\usepackage{amssymb,amsfonts}
\usepackage{amscd}
\usepackage{amsfonts}
\usepackage{amsmath}
\usepackage{amssymb}
\usepackage{amsthm}
\usepackage{bbm}
\usepackage{CJK}
\usepackage{fancyhdr}
\usepackage{graphicx}
\usepackage{hyperref}
\usepackage{indentfirst}
\usepackage{latexsym}
\usepackage{mathrsfs}
\usepackage{xypic}

\usepackage[top=1in,bottom=1in,left=1.25in,right=1.25in]{geometry}
\def\<{\langle}
\def\>{\rangle}
\def\a{\alpha}
\def\b{\beta}
\def\ci{\circ}
\def\c{\cdot}
\def\d{\delta}
\def\D{\Delta}

\def\g{\gamma}

\def\i{\iota}

\def\r{\rho}
\def\lr{\longrightarrow}

\def\o{\otimes}
\def\om{\omega}

\def\th{\theta}

\def\s{\sigma}
\def\v{\varepsilon}
\def\vp{\varphi}

\def\<{\langle}
\def\>{\rangle}

\date{}
\begin{document}
\renewcommand{\baselinestretch}{1.2}
\renewcommand{\arraystretch}{1.0}
\title{\bf A Note on Braided $T$-categories over Monoidal
  Hom-Hopf Algebras}
\date{}
\author {{\bf Miman You \quad  Shuanhong Wang \footnote {Corresponding author:  Shuanhong Wang, shuanhwang2002@yahoo.com}}\\
{\small Department of Mathematics, Southeast University}\\
{\small Nanjing, Jiangsu 210096, P. R. of China}}
 \maketitle
\begin{center}
\begin{minipage}{12.cm}

 \noindent{\bf Abstract}  Let  ${\sl Aut}_{mHH}(H)$ denote the set
  of all automorphisms of a monoidal Hom-Hopf algebra $H$
  with bijective antipode in the sense of Caenepeel and Goyvaerts
  \cite{CG2011}. The main aim of this paper is to provide new examples
  of braided $T$-category in the sense of Turaev \cite{T2008}.
   For this, first we construct a monoidal Hom-Hopf
   $T$-coalgebra $\mathcal{MHD}(H)$
    and prove that the $T$-category $Rep(\mathcal{MHD}(H))$
  of representation of $\mathcal{MHD}(H)$
 is isomorphic to $\mathcal {MHYD}(H)$ as braided $T$-categories,
 if $H$ is finite-dimensional.
 Then we construct a new braided $T$-category
 $\mathcal{ZMHYD}(H)$ over $\mathbb{Z},$ generalizing the main
 construction by Staic \cite{S2007}.
\\

\end{minipage}
\end{center}
\section*{0. INTRODUCTION}

Braided $T$-categories introduced by Turaev \cite{T2008} are of interest due to their applications in homotopy quantum field theories,
 which are generalizations of ordinary topological quantum field theories. As such, they are interesting to different research communities
 in mathematical physics (see \cite{FY1989, K2004, T1994, VA2001, VA2005}).
  Although Yetter-Drinfeld modules over Hopf algebras provide examples of such braided
  $T$-categories, these are rather trivial.
  The wish to obtain more interesting homotopy quantum field theories provides a strong
  motivation to find new examples of braided $T$-categories.
\\

 The aim of this article is to construct new examples of braided $T$-categories isomorphic
 to the $T$-category $\mathcal {MHYD}(H)$ in \cite{YW2014}.
 For this purpose, we prove that, if $(H,A,H)$ is a Yetter-Drinfeld Hom-datum
 (the second $H$ is regarded as an $H$-Hom-bimodule coalgebra) in \cite{YW2014},
  with $H$ finite dimensional, then the category $ _{A}\mathcal{MHYD}^{H}(H)$
of Yetter-Drinfeld Hom-modules is isomorphic to the category of left modules
over the diagonal crossed Hom-product $H^{*}\bowtie A$. Then
when $H$ is finite-dimensional we construct a monoidal Hom-Hopf $T$-coalgebra
 $\mathcal{MHYD}(H)$, and prove that the $T$-category
 $Rep(\mathcal{MHD}(H))$ of representation of $\mathcal{MHD}(H)$
 is isomorphic to $\mathcal {MHYD}(H)$ as braided $T$-categories.
\\

The article is organized as follows.
\\

  We will present the background material in Section 1.
  This section contains the relevant definitions on braided
  $T$-categories, monoidal Hom-Hopf algebras and monoidal Hom-Hopf
  $T$-coalgebras necessary for the understanding of the construction.
  In Section 2, we define the notion of a diagonal crossed
Hom-product algebra over a monoidal Hom-Hopf algebra.
And then when $H$ is finite dimensional,
we prove the category $ _{A}\mathcal{MHYD}^{H}(H)$
is isomorphic to the category of left $H^{*}\bowtie A$-modules,
 $_{H^{*}\bowtie A}\mathcal{M}$.
\\

  Section 3, when $H$ is finite-dimensional we construct a monoidal Hom-Hopf $T$-coalgebra
 $\mathcal{MHYD}(H)$, and prove that the $T$-category
 $Rep(\mathcal{MHD}(H))$ of representation of $\mathcal{MHD}(H)$
 is isomorphic to $\mathcal {MHYD}(H)$ as braided $T$-categories.

 Section 4, we construct a new braided $T$-category
 $\mathcal{ZMHYD}(H)$ over $\mathbb{Z},$ generalizing the main
 construction by Staic \cite{S2007}.

\section*{1. PRELIMINARIES}
\def\theequation{1. \arabic{equation}}
\setcounter{equation} {0} \hskip\parindent

Throughout, let $k$ be a fixed field. Everything is over $k$ unless
 otherwise specified.  We
 refer the readers to the books of Sweedler \cite{S1969}
  for the relevant concepts on the general theory of Hopf
 algebras.  Let $(C, \Delta )$ be a coalgebra. We use the "sigma" notation for
 $\Delta $ as follows:
 $$
 \Delta (c)=\sum c_1\otimes c_2, \,\,\forall c\in C.
 $$

\vskip 0.5cm
 {\bf 1.1. Braided $T$-categories.}
\vskip 0.5cm

 A {\sl monoidal category} ${\cal C}=({\cal C},\mathbb{I},\otimes,a,l,r)$
 is a category ${\cal C}$ endowed with a functor
 $\otimes: {\cal C}\times{\cal C}\rightarrow{\cal C}$
 (the {\sl tensor product}), an object $\mathbb{I}\in {\cal C}$
 (the {\sl tensor unit}), and natural isomorphisms $a$
 (the {\sl associativity constraint}), where
 $a_{U,V,W}:(U\otimes V)\otimes W\rightarrow U\otimes (V\otimes W)$
 for all $U,V,W\in {\cal C}$, and $l$ (the {\sl left unit constraint})
 where $l_U: \mathbb{I}\otimes U\rightarrow U,\,r$
 (the {\sl right unit constraint}) where
 $r_{U}:U\otimes{\cal C}\rightarrow U$ for all $U\in {\cal C}$,
 such that for all $U,V,W,X\in {\cal C},$
 the {\sl associativity pentagon}
 $a_{U,V,W\otimes X}\circ a_{U\otimes V,W,X}
 =(U\otimes a_{V,W,X})\circ a_{U,V\otimes W,X}\circ
 (a_{U,V,W}\otimes X)$ and
 $(U\otimes l_V)\circ(r_U\otimes V)=a_{U,I,V}$ are satisfied.
 A monoidal categoey ${\cal C}$ is {\sl strict} when all
  the constraints are identities.
\\

  Let $G$ be a group and let $Aut({\cal C})$ be the group of
  invertible strict tensor functors from ${\cal C}$ to itself.
   A category ${\cal C}$ over $G$
  is called a {\sl crossed category } if it satisfies the following:
  \begin{eqnarray*}
  &\blacklozenge & {\cal C} \mbox{ is a monoidal category;}\\
  &\blacklozenge & {\cal C} \mbox{ is disjoint union of a family
  of subcategories }\{{\cal C}_{\a }\}_{\a \in
  G},\mbox{ and for any }U\in {\cal C}_{\a },\\
  &&V\in {\cal C}_{\b }, U\o V\in {\cal C}_{\a \b }.
  \mbox{ The subcategory }{\cal C}_{\a }
  \mbox{ is called the }\a\mbox{th component of }{\cal C};\\
  &\blacklozenge & \mbox{Consider a group homomorphism }
   \vp : G\lr Aut({\cal C}), \b \mapsto \vp _{\b }, \mbox{ and
  assume that}\\
  &&  \vp _{\b }(\vp _{\a })
  =\vp _{\b\alpha\beta^{-1}},
  \mbox{ for all }\alpha,\beta\in G.\mbox{ The functors }
  \vp _{\b } \mbox{ are called conjugation}\\
  &&\mbox{ isomorphisms.}
  \end{eqnarray*}

 Furthermore, $ {\cal C}$ is called strict when it is
 strict as a monoidal category.
\\

 {\sl Left index notation}: Given $\a \in G$
 and an object $V\in {\cal C}_{\a }$, the functor $\vp _{\a }$
 will be denoted by ${}^V( \cdot )$, as in Turaev \cite{T2008} or
 Zunino \cite{Z2004}, or even ${}^{\a }( \cdot )$.
 We use the notation ${}^{\overline{V}}( \cdot )$
 for ${}^{\a ^{-1}}( \cdot )$. Then we have
 ${}^V id_U=id_{{} V^U}$ and
 ${}^V(g\circ f)={}^Vg\circ {}^Vf$.
 Since the conjugation $\vp : G\lr Aut({\cal C})$ is a
 group homomorphism, for all $V, W\in {\cal C}$, we have ${}^{V\o W}( \cdot )
 ={}^V({}^W( \cdot ))$ and ${}^\mathbb{I}( \cdot )={}^V({}^{\overline{V}}( \cdot ))
 ={}^{\overline{V}}({}^V( \cdot ))=id_{\cal C}$. Since, for all
 $V\in {\cal C}$, the functor ${}^V( \cdot )$ is strict, we have
 ${}^V(f\o g)={}^Vf\o {}^Vg$, for any morphisms $f$ and $g$ in ${\cal C}$,
  and ${}^V\mathbb{I}=\mathbb{I}$.
\\

 A {\sl braiding} of a crossed category ${\cal C}$ is
 a family of isomorphisms $({c=c_{U,V}})_{U,V}\in {\cal C}$,
 where $c_{U,V}: U\otimes V\rightarrow {}^UV\otimes U$
 satisfying the following conditions:
 \begin{itemize}
 \item [(1)] For any arrow $f\in {\cal C}_{\a }(U, U')$ and
  $g\in {\cal C}(V, V')$,
  $$
 (({}^{\a }g)\o f)\circ c _{U, V}=c _{U' V'}\circ (f\o g).
  $$
 \item [(2)]
  For all $ U, V, W\in {\cal C},$ we have
  \begin{eqnarray*}
  && c _{U\o V, W}=a_{{}^{U\o V}W, U, V}\circ (c _{U, {}^VW}\o
  id_V)\circ a^{-1}_{U, {}^VW, V}\circ (\i _U\o c _{V, W})
  \circ a_{U, V, W},\\
  && c _{U, V\o W}=a^{-1}_{{}^UV, {}^UW, U}
 \circ (\i _{({}^UV)}\o c _{U, W})\circ a_{{}^UV, U, W}\ci
 (c _{U, V}\o \i_W)\circ a^{-1}_{U, V, W},
 \end{eqnarray*}
 where $a$ is the natural isomorphisms in the tensor category
 ${\cal C}$.
 \item [(3)] For all $ U, V\in {\cal C}$ and $\b\in G$,
 $$ \vp _{\b }(c_{U, V})=c _{\vp _{\b }(U), \vp _{\b }(V)}. $$
 \end{itemize}

 A crossed category endowed with a braiding is called
 a {\sl braided $T$-category}.
\\

{\bf 1.2. Monoidal Hom-Hopf algebras.}
\vskip 0.5cm

Let $\mathcal{M}_{k}=(\mathcal{M}_{k},\o,k,a,l,r )$
 denote the usual monoidal category of $k$-vector spaces and linear maps between them.
 Recall from \cite{CG2011}
 that there is the {\it monoidal Hom-category} $\widetilde{\mathcal{H}}(\mathcal{M}_{k})=
 (\mathcal{H}(\mathcal{M}_{k}),\,\o,\,(k,\,id),
 \,\widetilde{a},\,\widetilde{l},\,\widetilde{r })$, a new monoidal category,
  associated with $\mathcal {M}_{k}$ as follows:
\begin{itemize}
 \item The objects of the monoidal category
 $ \mathcal{H}(\mathcal{M}_{k})$ are couples
 $(M,\xi_{M})$, where $M \in \mathcal {M}_{k}$ and $\xi_{M} \in Aut_k(M)$, the set of
   all $k$-linear automomorphisms of $M$;

  \item The morphism $f:(M,\xi_{M})\rightarrow (N,\xi_{N})$ in $ \mathcal{H}(\mathcal{M}_{k})$
  is the $k$-linear map $f: M\rightarrow N$ in $\mathcal{M}_{k}$
  satisfying   $ \xi_{N} \circ f = f\ci \xi_{M}$, for any two objects
  $(M,\xi_{M}),(N,\xi_{N})\in \mathcal{H}(\mathcal{M}_{k})$;

 \item The tensor product is given by
 $$
 (M,\xi_{M})\o (N,\xi_{N})=(M\o N,\xi_{M}\o\xi_{N} )
 $$
for any $(M,\xi_{M}),(N,\xi_{N})\in \mathcal{H}(\mathcal{M}_{k})$.

\item The tensor unit is given by $(k, id)$;

\item The associativity constraint $\widetilde{a}$
 is given by the formula
 $$
 \widetilde{a}_{M,N,L}=a_{M,N,L}\circ((\xi_{M}\o id)\o
  \xi_{L}^{-1})=(\xi_{M}\o(id\o\xi_{L}^{-1}))\circ a_{M,N,L},
  $$
 for any objects
 $(M,\xi_{M}),(N,\xi_{N}),(L,\xi_{L})\in \widetilde{\mathcal{H}}(\mathcal{M}_{k})$;

 \item The left and right unit constraint
  $\widetilde{l}$ and $\widetilde{r }$ are given by
 $$
 \widetilde{l}_M=\xi_{M}\circ l_M=l_M\circ(id\o\xi_{M}),\, \quad
  \widetilde{r}_M =\xi_{M}\circ r_M=r_M\circ(\xi_{M}\o id)
  $$
for all $(M,\xi_{M}) \in \widetilde{\mathcal{H}}(\mathcal{M}_{k})$.
\end{itemize}

We  now recall from \cite{CG2011} the following notions used later.
\\

  {\it Definition 1.2.1.}  Let $\widetilde{\mathcal{H}}(\mathcal{M}_{k})$
  be a monoidal Hom-category. A {\it monoidal Hom-algebra} is an object
  $(A,\xi_{A})$ in $\widetilde{\mathcal{H}}(\mathcal{M}_{k})$ together
 with an element $1_A\in A$ and linear maps
 $$m:A\o A\rightarrow A;\,\,a\o b\mapsto ab, \,\,\,\xi_{A}\in Aut_k(A)$$
 such that
  \begin{eqnarray}
 \xi_{A}(ab)=\xi_{A}(a)\xi_{A}(b),&& \alpha(1_A)=1_A ,\label{eq1}\\
  \xi_{A}(a)(bc)=(ab)\xi_{A}(c),&& a1_A=1_Aa=\xi_{A}(a),\label{eq2}
  \end{eqnarray}
for all $a,b,c\in A.$
\\

  {\it Definition 1.2.2. A monoidal Hom-coalgebra} is
  an object $(C,\xi_{C})$ in the category
  $\widetilde{\mathcal{H}}(\mathcal{M}_{k})$
 together with linear maps
 $\D:C\rightarrow C\o C,\,\D(c)=c_1\o c_2$ and
 $\varepsilon:C\rightarrow k$ such that
 \begin{eqnarray}
  \D(\xi_{C}(c))=\xi_{C}(c_1)\o\xi_{C}(c_2),
  && \varepsilon(\xi_{C}(c))=\varepsilon(c),\label{eq3}\\
   \xi_{C}^{-1}(c_1)\o\D(c_2)=\D(c_1)\o\xi_{C}^{-1}(c_2),
   && c_1\varepsilon(c_2)=\xi_{C}^{-1}(c)=\varepsilon(c_1)c_2,\label{eq4}
  \end{eqnarray}
 for all $c\in C.$
\\

 {\bf Remark 1.2.3.} (1) Note that (1.4) is equivalent to
 $c_1\o c_{21}\o \xi_{C}(c_{22})=\xi_{C}(c_{11})\o c_{12}\o c_2.$
 Analogue to monoidal Hom-algebras, monoidal Hom-coalgebras
 will be short for counital monoidal Hom-coassociative coalgebras
 without any confusion.

 (2)  Let $(C,\xi_{C})$ and $(C',\xi_{C}')$ be two monoidal Hom-coalgebras.
 A monoidal Hom-coalgebra map $f:(C,\xi_{C})\rightarrow(C',\xi_{C}')$
 is a linear map such that $f\circ \xi_{C}=\xi_{C}'\circ f, \D\circ f=(f\o f)\circ\D$
  and $\varepsilon'\circ f=\varepsilon.$
\\

 {\it Definition 1.2.4. A monoidal Hom-Hopf algebra} $H=(H,\xi_{H},m,1_H,\D,\varepsilon,S)$
 is a bialgebra with $S$ in $ \widetilde{\mathcal{H}}(\mathcal {M}_{k}).$
 This means that $(H,\alpha,m,1_H)$ is a monoidal Hom-algebra and
 $(H,\alpha,\D,\varepsilon)$ is
 a monoidal Hom-coalgebra such that $\D$ and $\varepsilon$
  are morphisms of algebras,
 that is, for all $h,g\in H,$
 $$\D(hg)=\D(h)\D(g),\, \,\, \D(1_H)=1_H\o1_H,\,\,\,\,\,\,
  \varepsilon(hg)=\varepsilon(h)\varepsilon(g), \,\,\,\,\,\varepsilon(1_H)=1.$$
$S$ is the convolution inverse of the identity morphism $id_H$
 (i.e., $ S*id=1_H\ci \varepsilon=id*S$). Explicitly,  for all $h\in H$,
 $$
 S(h_1)h_2=\varepsilon(h)1_H=h_1S(h_2).
 $$

{\bf Remark 1.2.5.} (1) Note that a monoidal Hom-Hopf algebra is
 by definition a Hopf algebra in $ \widetilde{\mathcal{H}}(\mathcal {M}_{k})$.

 (2)  Furthermore, the antipode of monoidal Hom-Hopf algebras has
   almost all the properties of antipode of Hopf algebras such as
 $$S(hg)=S(g)S(h),\,\,\,\, S(1_H)=1_H,\,\,\,\,
  \D(S(h))=S(h_2)\o S(h_1),\,\,\,\,\,\,\varepsilon\ci S=\varepsilon.$$
 That is, $S$ is a monoidal Hom-anti-(co)algebra homomorphism.
  Since $\xi_{H}$ is bijective and commutes with $S$,
  we can also have that the inverse $\xi_{H}^{-1}$ commutes with $S$,
  that is, $S\ci \xi_{H}^{-1}= \xi_{H}^{-1}\ci S.$
 \\

 In the following, we recall the notions of actions on monoidal
 Hom-algebras and coactions on monoidal Hom-coalgebras.
\\

  Let $(A,\xi_{A})$ be a monoidal Hom-algebra.
{\it A left $(A,\xi_{A})$-Hom-module} consists of
 an object $(M,\xi_{M})$ in $\tilde{\mathcal{H}}(\mathcal {M}_{k})$
  together with a morphism
  $\psi:A\o M\rightarrow M,\psi(a\o m)=a\cdot m$ such that
 $$\xi_{A}(a)\c(b\c m)=(ab)\c\xi_{M}(m),\,\,
 \,\,\xi_{M}(a\c m)=\xi_{A}(a)\c\xi_{M}(m),\,\,
 \,\,1_A\c m=\xi_{M}(m),$$
 for all $a,b\in A$ and $m \in M.$
\\

 Monoidal Hom-algebra $(A,\xi_{A})$ can be
 considered as a Hom-module on itself by the Hom-multiplication.
  Let $(M,\xi_{M})$ and $(N,\xi_{N})$ be two left $(A,\xi_{A})$-Hom-modules.
  A morphism $f:M\rightarrow N$ is called left
   $(A,\xi_{A})$-linear if
   $f(a\c m)=a\c f(m),f\ci \xi_{M}= \xi_{N}\ci f.$
  We denoted the category of left $(A,\xi_{A})$-Hom modules by
  $\widetilde{\mathcal{H}}(_{A}\mathcal {M}_{k})$.
\\

 Similarly, let $(C,\xi_{C})$ be a monoidal Hom-coalgebra.
 {\it A right $(C,\xi_{C})$-Hom-comodule} is an object
  $(M,\xi_{M})$ in $\tilde{\mathcal{H}}(\mathcal {M}_{k})$
  together with a $k$-linear map
  $\rho_M:M\rightarrow M\o C,\rho_M(m)=m_{(0)}\o m_{(1)}$ such that
 \begin{equation}
 \xi_{M}^{-1}(m_{(0)})\o \D_C(m_{(1)})
 =(m_{(0)(0)}\o m_{(0)(1)})\o \xi_{C}^{-1}(m_{(1)}),\label{eq5}
  \end{equation}
  \begin{equation}\label{eq6}
    \rho_M(\xi_{M}(m))=\xi_{M}(m_{(0)})\o\xi_{C}(m_{(1)}),
     \ \ \
   m_{(0)}\varepsilon(m_{(1)})=\xi_{M}^{-1}(m),
  \end{equation}
for all $m\in M.$
\\

 $(C,\xi_{C})$ is a Hom-comodule on itself via the Hom-comultiplication.
 Let $(M,\xi_{M})$ and $(N,\xi_{N})$ be two right $(C,\xi_{C})$-Hom-comodules.
  A morphism $g:M\rightarrow N$ is called right $(C,\xi_{C})$-colinear
  if $g\ci \xi_{M}=\xi_{N}\ci g$ and
  $g(m_{(0)})\o m_{(1)}=g(m)_{(0)}\o g(m)_{(1)}.$
  The category of right
  $(C,\xi_{C})$-Hom-comodules is denoted by
   $\widetilde{\cal{H}}(\cal {M}^C)$ .
\\

 {\it Definition 1.2.6.} Let $(H,m,\Delta,S,\xi_{H})$ be a monoidal Hom-bialgebra
 and $\a, \b\in  {\sl Aut}_{mHH}(H)$. Recall from \cite{YW2014}
 that a {\it left-right $(\a, \b)$-Yetter-Drinfeld Hom-module } over $(H,\xi_{H})$
 is the object $(M,\cdot,\rho,\xi_{M})$ which
 is both in  $\widetilde{\cal{H}}(_{H}\cal {M})$ and $\widetilde{\cal{H}}(\cal {M}^{H})$
 obeying the compatibility condition:
 \begin{equation}\label{eq7}
\r (h\c m)=\xi_{H}(h_{21})\c m_{0}\o (\b(h_{22})\xi_{H}^{-1}(m_{1}))\a(S^{-1}(h_{1})),
\end{equation}

 {\bf Remark 1.2.7.} (1) The category of all left-right $(\a,\b)$-Yetter-Drinfeld Hom-modules
  is denoted by $ _{H}\mathcal{MHYD}^{H}(\a, \b)$ with understanding morphism.

 (2) If $(H,\xi_{H})$ is a monoidal Hom-Hopf algebra with a bijective
 antipode $S$ and $S$ commute with $\a,\b$, then the above equality is equivalent to
 \begin{equation}\label{eq8}
h_{1}\c m_{0}\o \b(h_{2})m_{1}=\xi_{M}((h_{2}\c\xi_{M}^{-1} (m))_{0}) \o (h_{2}\c
 \xi_{M}^{-1}(m))_{1}\a(h_{1}).
\end{equation}
 for all $h\in H$ and $m\in M$.

 (3) If $(M,\xi_{M}) \in {_{H}}\mathcal {MHYD}^{H}(\a,\b)$
 and $(N,\xi_{N})\in {_{H}}\mathcal {MHYD}^{H}(\gamma,\delta)$, with
 $\a,\b,\g,\d \in {\sl Aut}_{mHH}(H)$, then $(M \o N,\xi_{M}\o\xi_{N})
 \in {_{H}}\mathcal {MHYD}^{H}(\a\g, \d\g^{-1}\b\g)$ with structures as follows:
\begin{eqnarray}
h\c (m \o n) &=& \g (h_{1})\c m \o \g^{-1}\b\g(h_{2})\c n,\label{eq12}\\
m\o n &\mapsto & (m_{0}\o n_{0})\o n_{1}m_{1}.\label{eq13}
\end{eqnarray}
for all $m\in M,n\in N$ and $h\in H.$
\\

  {\it Definition 1.2.8.} Let $(H,\xi_{H})$ be a monoidal Hom-algebra. A monoidal Hom-algebra
  $(A,\xi_{A})$ is called an {\it $(H,\xi_{H})$-Hom-bicomodule
  algebra} in \cite{YW2014}, with Hom-comodule maps $\rho_l$ and $\rho_r$
  obeying the following axioms:
 \begin{itemize}
  \item [(1)] $\rho_l: A\rightarrow H\o A,\ \rho_l(a) = a_{[-1]}\otimes a_{[0]}$, and
  $\rho_r: A\rightarrow A\o H,\ \rho_r(a)= a_{<0>}\otimes a_{<1>},$
  \item [(2)] $\rho_l$ and $\rho_r$ satisfy the following compatibility condition:
  for all $a\in A,$
  \begin{eqnarray}\label{eq9}
  a_{<0>[-1]}\otimes a_{<0>[0]}\otimes \xi_{H}^{-1}(a_{<1>})
  = \xi_{H}^{-1}(a_{[-1]})\otimes a_{[0]<0>}\otimes a_{[0]<1>}.
    \end{eqnarray}
 \end{itemize}

  {\it Definition 1.2.9.} Let $(H,\xi_{H})$ be a monoidal Hom-Hopf algebra,
 $(A,\xi_{A})$ be an $H$-Hom-bicomodule algebra.
 We consider {\it the Yetter-Drinfeld Hom-datum $(H,A,H)$} as in \cite{YW2014},
 (the second $H$ is regarded as an $H$-Hom-bimodule coalgebra),
 and {\it the Yetter-Drinfeld Hom-module category $ _{A}\mathcal{MHYD}^{H}(H)$,}
 whose objects are $k$-modules $(M,\xi_{M})$ with the following additional structure:
 \begin{itemize}
  \item [(1)] $M$ is a left $A$-module;
  \item [(2)] we have a $k$-linear map $\rho_{M}:M\rightarrow M\o H, \rho_{M}(m)=m_{0}\otimes m_{1},$
  \item [(3)] the following compatibility conditions holds:
  \begin{eqnarray}
 (a\cdot m)_{0}\otimes(a\cdot m)_{1}
  =\xi_{A}(a_{[0]<0>})\cdot m_{0}
 \otimes( a_{[0]<1>}\alpha^{-1}(m_{1}))S^{-1}(a_{[-1]}),&& \label{eq10}\\
  a_{<0>}\cdot m_{0}\otimes a_{<1>} m_{1}
  =\xi_{M}((a_{[0]}\cdot \xi_{M}^{-1}(m))_{0})
 \otimes(a_{[0]}\cdot \xi_{M}^{-1}(m))_{1}a_{[-1]}. &&\label{eq11}\ \
\end{eqnarray}
 \end{itemize}
  for all $a\in A$ and $m\in M$. \\

  {\it Definition 1.2.10.} Let $(A, \xi_A)$ be a monoidal Hom-algebra,
  $(M,\xi_{M})$ be a monoidal Hom-algebra. Assume that $(M, \alpha _M)$ is both a left and a right
$A$-module algebra (with actions denoted by $A\o M\rightarrow M$, $a\o m\mapsto a\cdot m$ and $M\o A
\rightarrow M$, $m\o a\mapsto m\cdot a$). We call $(M, \xi _M)$ an {\em $A$-bimodule} as in
\cite{FP2014} if the following condition is satisfied, for all $a, a'\in A$, $m\in M$:
\begin{eqnarray}\label{bimodule}
&&\xi _A(a)\cdot (m\cdot a')=(a\cdot m)\cdot \xi _A(a').
\label{hombimodule}
\end{eqnarray}

{\bf 1.3. Monoidal Hom-Hopf $T$-coalgebras.}
\vskip 0.5cm

 {\it Definition 1.3.1.} Let $G$ be a group with unit $1$. Then we recall from Yang Tao \cite{Y2014}
 that {\it a monoidal Hom $T$-coalgebra}
  $(C, \xi_{C})$ over $G$ is a family of objects $\{(C_{p}, \xi_{C_{p}})\}_{p\in G}$
 in $\mathcal{\widetilde{H}}(\mathcal{M}_{k})$ together with linear maps
 $\Delta_{p, q}: C_{pq}\longrightarrow C_{p}\otimes C_{q}, c_{pq}\mapsto c_{(1, p)}\otimes c_{(2, q)}$
 and $\varepsilon: C_{e}\longrightarrow k$ such that
 \begin{eqnarray*}
 \xi_{C_{P}}^{-1}(c_{(1, p)}) \otimes \Delta_{q, r}(c_{(2, qr)}) = \Delta_{p, q}(c_{(1, pq)}) \otimes \xi_{c_{r}}^{-1}(c_{(2, r)}),
 && \forall c\in C_{pqr}, \\
  c_{(1, p)}\varepsilon(c_{(2, e)})=\varepsilon(c_{(1, e)})c_{(2, p)}=\xi_{C_{p}}^{-1}(c_{p}), && \forall c\in C_{p}, \\
 \Delta_{p, q}(\xi{C_{pq}}^{-1}(c_{pq}))=\xi_{C_{p}}^{-1}(c_{(1, p)}) \otimes \xi_{C_{q}}^{-1}(c_{(2, q)}), && \forall c\in C_{pq}, \\
 \varepsilon(\xi_{C_{e}}^{-1}(c))=\varepsilon(c), && \forall c\in C_{e}.
 \end{eqnarray*}

 Let $(C, \xi_{C})$ and $(C', \xi_{C}')$ be two monoidal Hom $T$-coalgebras over $G$.
 A Hom-coalgebra map $f: (C, \xi_{C})\longrightarrow(C', \xi_{C}')$
 is a family of linear maps $\{f_{p}\}_{p\in G}$, $f_{p}: (C_{p}, \xi_{C_{p}})\longrightarrow(C'_{p}, \xi_{C_{p}}')$
 such that $f_{p}\circ \xi_{C_{p}} = \xi_{C_{p}}'\circ f_{p}$, $\Delta_{p, q}\circ f_{pq} = (f_{p}\otimes f_{q})\Delta_{p, q}$ and $\varepsilon\circ f_{e}= \varepsilon$.
 \\

  {\it Definition 1.3.2.} A {\it monoidal Hom-Hopf $T$-coalgebra} $(H=\bigoplus_{p\in G}H_{p},
 \xi= \{\xi_{H_{p}}\}_{p\in G})$
 is a monoidal Hom $T$-coalgebra
 where each $(H_{p}, \xi_{H_{p}})$ is a monoidal Hom-algebra with multiplication $m_{p}$ and unit $1_{p}$
 endowed with antipode $S=\{S_{p}\}_{p\in G}$,
 $S_{p}: H_{p}\longrightarrow H_{p^{-1}} \in \mathcal{\widetilde{H}}(\mathcal{M}_{k})$ such that
 \begin{eqnarray*}
 \Delta_{p, q}(hg) = \Delta_{p, q}(h)\Delta_{p, q}(g), \quad \Delta_{p, q}(1_{pq}) = 1_{p} \otimes 1_{q}, && \forall h, g\in H_{pq}\\
 \varepsilon(hg)=\varepsilon(h)\varepsilon(g), \qquad \varepsilon(1_{e})=1_{k}, && \forall h, g\in H_{e}\\
 S_{p^{-1}}(h_{(1, p^{-1})})h_{(2, p)} = \varepsilon(h)1_{p} = h_{(1, p)}S_{p^{-1}}(h_{(2, p^{-1})}) && \forall h \in H_{e}.
 \end{eqnarray*}

 Note also that the $(H_{e}, \xi_{e}, m_{e}, 1_{e}, \Delta_{e, e}, \varepsilon, S_{e})$
  is a monoidal Hom-Hopf algebra in the usual sense of the word.
   We call it the neutral component of $H$.

 {\it Definition 1.3.3.} A monoidal Hom-Hopf $T$-coalgebra $(H=\bigoplus_{p\in G}H_{p},
 \xi= \{\xi_{H_{p}}\}_{p\in G})$ is called a {\it monoidal Hom-Hopf crossed monoidal Hom-Hopf
 $T$-coalgebra } if it is endowed with a family of
 algebra isomorphisms $\varphi= \{\varphi_{\beta}^{\alpha}:
  H_{\alpha}\rightarrow H_{\beta\alpha\beta^{-1}}\}_{\alpha,\beta\in G}$
  such that
\begin{itemize}
    \item each $\varphi_{\g}$ preserves the comultiplication and the counit i.e.,
    for any $\alpha,\beta,\gamma\in G,$
    we have $\Delta_{\gamma\alpha\gamma^{-1},\gamma\beta\gamma^{-1}}
   \circ\varphi_{\gamma}=(\varphi_{\gamma}\otimes\varphi_{\gamma})
   \circ\Delta_{\alpha,\beta}$ and $\varepsilon \circ\varphi_{\gamma}=\varepsilon.$

   \item $\varphi$ is multiplicative, i.e.,
   $\varphi_{\beta}\circ\varphi_{\gamma}=\varphi_{\beta\gamma},$
   for any $\beta,\gamma\in G.$
\end{itemize}

 It is easy to get the following identities,
 $\varphi_1|H_{\alpha}=id_{\alpha}$ and  $\varphi_{\alpha}^{-1}=\varphi_{\alpha^{-1}},$
 for all $\alpha\in G.$
 Moreover, $\varphi$ preserves the antipode, i.e.,
 $\varphi_{\beta}\circ S_{\alpha}
 =S_{\beta\alpha\beta^{-1}}\circ\varphi_{\beta}$
 for all $\alpha,\beta\in G$.
 \\

\section*{2. THE DIAGONAL CROSSED HOM-PRODUCT}
\def\theequation{2. \arabic{equation}}
\setcounter{equation}{0} \hskip\parindent

In this section, we define the notion of the diagonal crossed
Hom-product over a monoidal Hom-Hopf algebra that
are based on Hom-associative left and right coactions.
If $H$ is finite dimensional, we prove the category $ _{A}\mathcal{MHYD}^{H}(H)$
is isomorphic to the category of left $H^{*}\bowtie A$-modules,
 $_{H^{*}\bowtie A}\mathcal{M}$, generalizing the results in \cite{DPV2006}.
\\

In what follows, let  $(H,\xi_{H})$ be a monoidal Hom-Hopf algebra with
the bijective antipode $S$  and let ${\sl Aut}_{mHH}(H)$ denote the set
  of all automorphisms of a monoidal Hopf algebra $H$.
\\

{\bf Definition 2.1.} Let $(H,\xi_{H})$ be a finite dimensional monoidal Hom-Hopf algebra,
 $(A,\xi_{A})$ be a monoidal Hom-bicomodule algebra. Then the diagonal crossed Hom-product
  $H^{*}\bowtie A$ is defined as follows:
  \begin{itemize}
  \item [--] as $k$-spaces, $H^{*}\bowtie A=H^{*}\o A$;
  \item [--] multiplication is given by
  \begin{eqnarray}\label{2.2}
  (f\bowtie a)(g\bowtie b) = f(a_{[-1]}\rightharpoonup
  (\xi_{H}^{*2}(g)\leftharpoonup S^{-1}(a_{[0]_{<1>}})))\bowtie \xi_{A}^{2}(a_{[0]_{<0>}})b;
  \end{eqnarray}
  \begin{eqnarray}\label{2.1}
  h\rightharpoonup f = \<f_{2},\xi_{H}^{-1}(h)\>\xi_{H}^{*-2}(f_{1})\ \
   \mbox{ and } \ \ f\leftharpoonup h = \<f_{1},\xi_{H}^{-1}(h)\>\xi_{H}^{*-2}(f_{2});
  \end{eqnarray}
    \end{itemize}
 for all $a,b\in (A,\xi_{A}),\ f,g\in (H^{*},\xi_{H}^{*-1}),h\in (H,\xi_{H}).$
 \\

  {\bf Proposition 2.2.} Let $(A,\xi_{A})$ be an $(H,\xi_{H})$-Hom-bicomodule algebra
  and $(H^{*},\xi_{H}^{*-1})$ be an $(H,\xi_{H})$-Hom-bimodule algebra. Then the tensor space
  $H^{*}\o A$ is a Hom-algebra with the multiplication in the formula (\ref{2.2})
  and the unit $\varepsilon_{H}\bowtie 1_{A}.$
\\

{\bf Proof.} It is obvious that
$(\varepsilon_{H}\bowtie 1_{A})(f\bowtie a)=\xi_{H}^{*-1}(f)\bowtie \xi_{A}(a),$
so $(\varepsilon_{H}\bowtie 1_{A})$ is unit element. We have:

 \begin{eqnarray*}
  &&[(f\bowtie a)(g\bowtie b)]\xi_{H^*\bowtie A}(\phi\bowtie c)\\
  &=&[f(a_{[-1]}\rightharpoonup(\xi_{H}^{*2}(g)\leftharpoonup S^{-1}(a_{[0]_{<1>}})))]
  ((\xi_{A}^{2}(a_{[0]_{<0>}})b)_{[-1]}\rightharpoonup(\xi_{H}^{*}(\phi)\\
  &&\leftharpoonup S^{-1}((\xi_{A}^{2}(a_{[0]_{<0>}})b)_{[0]_{<1>}})))
  \bowtie\xi_{A}^{2}((\xi_{A}^{2}(a_{[0]_{<0>}})b)_{[0]_{<0>}})\xi_{A}(c)\\
  &=&\xi_{H}^{*-1}(f)[(a_{[-1]}\rightharpoonup(\xi_{H}^{*2}(g)\leftharpoonup S^{-1}(a_{[0]_{<1>}})))
  \xi_{H}^{*}((\xi_{H}^{2}(a_{[0]_{<0>_{[-1]}}})b_{[-1]})\rightharpoonup(\xi_{H}^{*}(\phi)\\
  &&\leftharpoonup S^{-1}(\xi_{H}^{2}(a_{[0]_{<0>_{[0]_{<1>}}}})b_{[0]_{<1>}})))]
  \bowtie\xi_{A}^{5}(a_{[0]_{<0>_{[0]_{<0>}}}})(\xi_{A}^{2}(b_{[0]_{<0>}})c)\\
  &=&\xi_{H}^{*-1}(f)[(a_{[-1]}\rightharpoonup(\xi_{H}^{*2}(g)\leftharpoonup S^{-1}(a_{[0]_{<1>}})))
    \<\xi_{H}^{*-1}(\phi_{221}),\xi_{H}(a_{[0]_{<0>_{[-1]}}})\> \\
 && \<\xi_{H}^{*-1}(\phi_{222}),\xi_{H}^{-1}(b_{[-1]})\>
    \<\xi_{H}^{*}(\phi_{11}),S^{-1}\xi_{H}^{-1}(b_{[0]_{<1>}})\>
    \<\xi_{H}^{*}(\phi_{12}),S^{-1}\xi_{H}(a_{[0]_{<0>_{[0]_{<1>}}}})\>\\
 && \xi_{H}^{*-2}(\phi_{21})]\bowtie\xi_{A}^{5}(a_{[0]_{<0>_{[0]_{<0>}}}})(\xi_{A}^{2}(b_{[0]_{<0>}})c)\\
 &=&\xi_{H}^{*-1}(f)[(a_{[-1]}\rightharpoonup(\xi_{H}^{*2}(g)\leftharpoonup S^{-1}(a_{[0]_{<1>}})))
    \<\xi_{H}^{*-2}(\phi_{2122}),\xi_{H}(a_{[0]_{<0>_{[-1]}}})\>\\
 && \<\phi_{22},\xi_{H}^{-1}(b_{[-1]})\>\<\xi_{H}^{*2}(\phi_{1}),S^{-1}\xi_{H}^{-1}(b_{[0]_{<1>}})\>
    \<\phi_{211},S^{-1}\xi_{H}(a_{[0]_{<0>_{[0]_{<1>}}}})\>\\
 &&  \xi_{H}^{*-4}(\phi_{2121})]\bowtie\xi_{A}^{5}(a_{[0]_{<0>_{[0]_{<0>}}}})
 (\xi_{A}^{2}(b_{[0]_{<0>}})c)\\
 &=&\xi_{H}^{*-1}(f)[(a_{[-1]}\rightharpoonup(\xi_{H}^{*2}(g)\leftharpoonup S^{-1}(a_{[0]_{<1>}})))
 (\xi_{H}^{2}(a_{[0]_{<0>_{[-1]}}})\rightharpoonup(\xi_{H}^{*2}(b_{[-1]}\rightharpoonup(\xi_{H}^{*2}(\phi)\\
 &&\leftharpoonup S^{-1}(b_{[0]_{<1>}})))\leftharpoonup S^{-1}\xi_{H}^{2}(a_{[0]_{<0>_{[0]_{<1>}}}})))]
  \bowtie\xi_{A}^{5}(a_{[0]_{<0>_{[0]_{<0>}}}})(\xi_{A}^{2}(b_{[0]_{<0>}})c)\\
 &=&\xi_{H}^{*-1}(f)[(\xi_{H}(a_{[-1]_{1}})\rightharpoonup(\xi_{H}^{*2}(g)\leftharpoonup S^{-1}\xi_{H}(a_{[0]_{<1>2}})))
 (\xi_{H}(a_{[-1]2})\rightharpoonup(\xi_{H}^{*2}(b_{[-1]}\rightharpoonup(\xi_{H}^{*2}(\phi)\\
 &&\leftharpoonup S^{-1}(b_{[0]_{<1>}})))\leftharpoonup S^{-1}\xi_{H}(a_{[0]_{<1>1}})))]
  \bowtie\xi_{A}^{3}(a_{[0]_{<0>}})(\xi_{A}^{2}(b_{[0]_{<0>}})c)\\
 &=&\xi_{H^*\bowtie A}(f\bowtie a)[(g\bowtie b)(\phi\bowtie c)].
  \end{eqnarray*}
  Thus the multiplication is Hom-associative. This completes the proof.\hfill $\blacksquare$
  \\

 {\bf Example 2.3.} (1) If $(A,\xi_{A})=(H,\xi_{H})$ and $\rho_l=\rho_r=\Delta$
 the formula (\ref{2.2}) coincides with the multiplication in the Drinfeld double $(D(H),\xi_{H}^{*-1}\o\xi_{H})=(H^{*cop}\bowtie H,\xi_{H}^{*-1}\o\xi_{H})$, i.e.
  \begin{eqnarray}\label{2.3}
  (f\bowtie h)(g\bowtie l) = f(h_{1}\rightharpoonup
  (\xi_{H}^{*2}(g)\leftharpoonup S^{-1}(h_{22})))\bowtie \xi_{A}^{2}(h_{21})l,
  \end{eqnarray}
   for all $f,g\in H^{*}$ and $h,l\in H.$
   \\

 (2) Recall from Example 2.5 in \cite{YW2014}
 that $\a,\b \in Aut_{mHH}(H)$ and
  as $k$-vector spaces $(H(\a,\b),\xi_{H})=(H,\xi_{H})$,
 and $(H(\a,\b),\xi_{H})\in\,\, _{H}\mathcal{MHYD}^{H}(\a, \b)$,
  with with right $H$-Hom-comodule structure via Hom-comultiplication
and left $H$-Hom-module structure given by:
$$
h\cdot x=(\b(h_2)\xi_{H}^{-1}(x))\a(S^{-1}(\xi_{H}(h_1))).
$$
for all $h,\,x\in H.$

The diagonal crossed product
$(A(\a,\b),\xi_{H}^{*-1}\o \xi_{H})=(H^{*}\bowtie H(\a,\b),\xi_{H}^{*-1}\o \xi_{H})$,
whose multiplication is
\begin{eqnarray}\label{2.10}
  (f\bowtie h)(g\bowtie l) = f(\a(h_{1})\rightharpoonup
  (\xi_{H}^{*2}(g)\leftharpoonup S^{-1}(\b(h_{22}))))\bowtie \xi_{H}^{2}(h_{21})l,
  \end{eqnarray}
  for all $f,g \in H^{*}$ and $h,l\in H$.
\\

   The Drinfeld double $D(H)$ is a Hom-Hopf algebra with coproduct
  $\Delta_{D(H)}$ given by
  \begin{eqnarray}\label{2.4}
  \Delta_{D(H)}(f\bowtie h) = (f_{2}\bowtie h_{1})\o(f_{1}\bowtie h_{2}),
  \end{eqnarray}
  for all $f\in H^{*}$ and $h\in H.$
 \\

{\bf Proposition 2.4.} Let $(A,\xi_{A})$ be an $(H,\xi_{H})$-Hom-bicomodule
algebra. Then $H^{*}\bowtie A$ is a $D(H)$-Hom-bicomodule algebra
with two coactions $\rho_{r_{D(H)}}: H^*\bowtie A\rightarrow (H^*\bowtie A)\o D(H)$
 and $\rho_{l_{D(H)}}:H^*\bowtie A\rightarrow  D(H)\o(H^*\bowtie A)$ given by
 \begin{eqnarray*}
  \rho_{r_{D(H)}}(f\bowtie a)&=& (f_{2}\bowtie a_{<0>})\o (f_{1}\o a_{<1>}),\\
   \rho_{l_{D(H)}}(f\bowtie a)&=& (f_{2}\bowtie a_{[-1]})\o (f_{1}\o a_{[0]}),
  \end{eqnarray*}
  where elements in $D(H)$ are written as $(f\o h),h\in H, f\in H^{*},a\in A$.
  \\

 {\bf Proof.} In view of (\ref{2.4}) the comodule axioms
 and the Hom-coassociative (\ref{eq9}) are obvious.
 We are left to prove that $\rho_{r_D(H)}$ and $\rho_{l_D(H)}$
  are Hom-algebra maps. To this end we use the following identities obviously
  holding for all $f\in H^{*},h,l\in H$
  \begin{eqnarray}\label{2.5}
  \rho(h\rightharpoonup(f\leftharpoonup l))
   = (\xi_{H}^{*-1}(f_{1})\leftharpoonup l)\o
   (\xi_{H}^{-1}(h)\rightharpoonup\xi_{H}^{*-1}(f_{2}) ),
  \end{eqnarray}
With this we now compute
 \begin{eqnarray*}
  &&\rho_{r_{D(H)}}(f\bowtie a)\rho_{r_D(H)}(g\bowtie b)\\
  &=&[(f_{2}\bowtie a_{<0>})\o(f_{1}\o a_{<1>})]
  [(g_{2}\bowtie b_{<0>})\o(g_{1}\o b_{<1>})]\\
  &=&(f_{2}(a_{<0>_{[-1]}}\rightharpoonup(\xi_{H}^{*2}(g_{2})\leftharpoonup
  S^{-1}(a_{<0>_{[0]_{<1>}}})))\bowtie\xi_{A}^{2}(a_{<0>_{[0]_{<0>}}})b_{<0>})\\
  &&\o (f_{1}(a_{<1>1}\rightharpoonup(\xi_{H}^{*2}(g_{1})\leftharpoonup
  S^{-1}(a_{<1>22})))\o\xi_{H}^{2}(a_{<1>21})b_{<1>})\\
  &=&(f_{2}\<g_{21},S^{-1}\xi_{H}(a_{<0>_{[0]_{<1>}}})\>
  \<g_{222},\xi_{H}^{-1}(a_{<0>_{[-1]}})\>\xi_{H}^{*-2}(g_{221})
  \bowtie\xi_{A}^{2}(a_{<0>_{[0]_{<0>}}})b_{<0>})\\
  && \o (f_{1}\<g_{11},S^{-1}\xi_{H}(a_{<1>22})\>\<g_{122},\xi_{H}^{-1}(a_{<1>1})\>
  \xi_{H}^{*-2}(g_{121})\o\xi_{H}^{2}(a_{<1>21})b_{<1>})\\
  &=&(f_{2}\<g_{22},a_{<0>_{[-1]}}\>\xi_{H}^{*-1}(g_{21})\bowtie\xi_{A}(a_{<0>_{[0]}})b_{<0>})
  \o (f_{1}\<g_{11},S^{-1}(a_{<1>2})\> \xi_{H}^{*-1}(g_{121})\\
  &&\o\xi_{H}(a_{<1>1})b_{<1>})\\
  &=&(f_{2}\<g_{22},\xi_{H}^{-1}(a_{[-1]})\>\xi_{H}^{*-1}(g_{21})
  \bowtie\xi_{A}^{2}(a_{[0]_{<0>_{<0>}}})b_{<0>})
  \o (f_{1}\<g_{11},S^{-1}(a_{[0]_{<1>}})\> \xi_{H}^{*-1}(g_{121})\\
  &&\o\xi_{H}^{2}(a_{[0]_{<0>_{<1>}}})b_{<1>})\\
  &=&(f_{2}(\xi_{H}^{-1}(a_{[-1]})\rightharpoonup \xi_{H}^{*}(g_{2}))
  \bowtie\xi_{A}^{2}(a_{[0]_{<0>_{<0>}}})b_{<0>})
  \o (f_{1}(\xi_{H}^{*}(g_{1})\leftharpoonup S^{-1}(a_{[0]_{<1>}}))\\
  &&\o\xi_{H}^{2}(a_{[0]_{<0>_{<1>}}})b_{<1>})\\
  &\stackrel{(\ref{2.5})}{=}&(f_{2}(a_{[-1]}\rightharpoonup (\xi_{H}^{*2}(g)
  \leftharpoonup S^{-1}(a_{[0]_{<1>}})))_{2}
  \bowtie\xi_{A}^{2}(a_{[0]_{<0>_{<0>}}})b_{<0>})\\
 && \o (f_{1}(a_{[-1]}\rightharpoonup (\xi_{H}^{*2}(g)
  \leftharpoonup S^{-1}(a_{[0]_{<1>}})))_{1}
  \o\xi_{H}^{2}(a_{[0]_{<0>_{<1>}}})b_{<1>})\\
  &=&\rho_{r_{D(H)}}(f(a_{[-1]}\rightharpoonup (\xi_{H}^{*2}(g)
  \leftharpoonup S^{-1}(a_{[0]_{<1>}})))
  \bowtie\xi_{A}^{2}(a_{[0]_{<0>}})b)\\
 &=& \rho_{r_{D(H)}}((f\bowtie a)(g\bowtie b)).
  \end{eqnarray*}
  Hence $\rho_{r_{D(H)}}$ is a Hom-algebra map. The argument for
  $\rho_{l_{D(H)}}$ is analogous.\hfill $\blacksquare$
\\

 {\bf Example 2.4.} Let $(H,\xi_{H})$ be finite dimensional.
 Then $(A(\a,\b),\xi_{H}^{*-1}\o\xi_{H})$ becomes a $D(H)$-bicomodule
 algebra, with structures
 \begin{eqnarray*}
  H^{*}\bowtie H(\a,\b)\rightarrow (H^*\bowtie H(\a,\b))\o D(H),
  && f\bowtie h\mapsto (f_{2}\bowtie h_{1})\o(f_{1}\bowtie \b(h_{2})),\\
   H^*\bowtie H(\a,\b)\rightarrow  D(H)\o(H^*\bowtie H(\a,\b)),
  &&f\bowtie h\mapsto (f_{2}\bowtie \a(h_{1}))\o(f_{1}\bowtie h_{2}).
  \end{eqnarray*}
for all $f\in H^{*},h\in H.$
\\

In the rest of this section we establish that if $(H,\xi_{H})$
is a monoidal Hom-Hopf algebra and is finite dimensional
then the category $ _{A}\mathcal{MHYD}^{H}(H)$ is isomorphic
to the category of left $H^{*}\bowtie A$-modules,
 $_{H^{*}\bowtie A}\mathcal{M}.$
\\

{\bf Lemma 2.5.} Let $(H,\xi_{H})$ be a monoidal Hom-Hopf algebra
and $(H,A,H)$ a Yetter-Drinfeld Hom-datum. We have a functor
$F: _{A}\mathcal{MHYD}^{H}(H)\rightarrow _{H^{*}\bowtie A}\mathcal{M},$
given by $F(M)=M$ as $k$-module, with the $H^{*}\bowtie A$-module structure
defined by
\begin{eqnarray}\label{2.6}
  (f\bowtie u)\triangleright m=\<f,(u\cdot\xi_{M}^{-1}(m))_{1}\>
  \xi_{M}^{2}((u\cdot\xi_{M}^{-1}(m))_{0}),
  \end{eqnarray}
  for all $f\in (H^{*},\xi_{H}^{*-1}),u\in (A,\xi_{A})$ and $m\in (M,\xi_{M})$.
  $F$ transforms a morphism to itself.
  \\

{\bf Proof.} For all $f,g\in H^{*},a,b\in A$ and $m\in M$, we compute:
 \begin{eqnarray*}
  &&[(f\bowtie a)(g\bowtie b)]\triangleright \xi_{M}(m)\\
  &=&[f(a_{[-1]}\rightharpoonup(\xi_{H}^{*2}(g)
  \leftharpoonup S^{-1}(a_{[0]_{<1>}})))\bowtie \xi_{A}^{2}(a_{[0]_{<0>}})b]
  \triangleright\xi_{M}(m)\\
  &=&\<g_{1},S^{-1}\xi_{H}(a_{[0]_{<1>}})\>\<g_{22},\xi_{H}^{-1}(a_{[-1]})\>
  [f\xi_{H}^{*2}(g_{21})\bowtie \xi_{A}^{2}(a_{[0]_{<0>}})b]\triangleright\xi_{M}(m)\\
  &=&\<g_{1},S^{-1}\xi_{H}(a_{[0]_{<1>}})\>\<g_{22},\xi_{H}^{-1}(a_{[-1]})\>
  \<f,((\xi_{A}^{2}(a_{[0]_{<0>}})b)\cdot m)_{11}\>\\
  && \<\xi_{H}^{*-2}(g_{21}),((\xi_{A}^{2}(a_{[0]_{<0>}})b)\cdot m)_{12}\>
  \xi_{M}^{2}(((\xi_{A}^{2}(a_{[0]_{<0>}})b)\cdot m)_{0})\\
   &\stackrel{( \ref{eq10})}{=}&
   \<g_{1},S^{-1}\xi_{H}(a_{[0]_{<1>}})\>\<g_{22},\xi_{H}^{-1}(a_{[-1]})\>
  \<f,((\xi_{H}^{2}(a_{[0]_{<0>_{[0]_{<1>}}}})b_{[0]_{<1>}})\xi_{H}^{-1}(m_{1}))_{1}\\
  &&S^{-1}(\xi_{H}^{2}(a_{[0]_{<0>_{[-1]}}})b_{[-1]})_{1}\>
  \<g_{21},\xi_{H}^{-2}(((\xi_{H}^{2}(a_{[0]_{<0>_{[0]_{<1>}}}})b_{[0]_{<1>}})\xi_{H}^{-1}(m_{1}))_{2}\\
 && S^{-1}(\xi_{H}^{2}(a_{[0]_{<0>_{[-1]}}})b_{[-1]})_{2})\>
  \xi_{M}^{2}(\xi_{A}(\xi_{A}^{2}(a_{[0]_{<0>_{[0]_{<0>}}}})b_{[0]_{<0>}})\cdot m_{0})\\
  &=&\<g,S^{-1}\xi_{H}(a_{[0]_{<1>}})
  ((((a_{[0]_{<0>_{[0]_{<1>2}}}}\xi_{H}^{-2}(b_{[0]_{<1>2}}))\xi_{H}^{-3}(m_{12}))\\
  &&S^{-1}(a_{[0]_{<0>_{[-1]1}}}\xi_{H}^{-2}(b_{[-1]1})))\xi_{H}^{-1}(a_{[-1]}))\>
  \<f,((\xi_{H}^{2}(a_{[0]_{<0>_{[0]_{<1>1}}}})b_{[0]_{<1>1}})\xi_{H}^{-1}(m_{11}))\\
  &&S^{-1}(\xi_{H}^{2}(a_{[0]_{<0>_{[-1]2}}})b_{[-1]2})\>
  \xi_{M}^{2}((\xi_{A}^{3}(a_{[0]_{<0>_{[0]_{<0>}}}})\xi_{A}(b_{[0]_{<0>}}))\cdot m_{0})\\
  &=&\<g,S^{-1}\xi_{H}^{3}(a_{[0]_{<1>22}})
  (((a_{[0]_{<1>21}}\xi_{H}^{-1}(b_{[0]_{<1>2}}))\xi_{H}^{-2}(m_{12}))(S^{-1}\xi_{H}^{-1}(b_{[-1]1})\\
  &&(S^{-1}\xi_{H}^{-1}(a_{[-1]12})\xi_{H}^{-1}(a_{[-1]11}))))\>
  \<f,((a_{[0]_{<1>1}})b_{[0]_{<1>1}})\xi_{H}^{-1}(m_{11}))\\
  &&S^{-1}(a_{[-1]2}b_{[-1]2})\>
  \xi_{M}^{2}(\xi_{A}(a_{[0]_{<0>}}b_{[0]_{<0>}}))\cdot m_{0})\\
  &=&\<g,\xi_{H}^{2}(b_{[0]_{<1>2}})(\xi_{H}^{-1}(m_{12})(S^{-1}(b_{[-1]1}))\>
  \<f,((\xi_{H}(a_{[0]_{<1>}})b_{[0]_{<1>1}})\xi_{H}^{-1}(m_{01}))\\
  &&S^{-1}(\xi_{H}^{-1}(a_{[-1]2})b_{[-1]2})\>
  \xi_{M}^{3}((a_{[0]_{<0>}}b_{[0]_{<0>}})\cdot m_{00})\\
  &=&\<g,(b_{[0]_{<1>}}\xi_{H}^{-2}(m_{1}))S^{-1}(b_{[-1]})\>
  \<f,((\xi_{H}^{-1}(a_{[0]_{<1>}})\xi_{H}(b_{[0]_{<0>_{[0]_{<1>}}}}))\xi_{H}^{-1}(m_{01}))\\
  &&S^{-1}(\xi_{H}^{-1}(a_{[-1]2})\xi_{H}(b_{[0]_{<0>_{[-1]}}}))\>
  \xi_{M}^{3}((a_{[0]_{<0>}}\xi_{A}(b_{[0]_{<0>_{[-1]}}}))\cdot m_{00})\\
  &\stackrel{(\ref{eq10})}{=}&
  \<g,(b_{[0]_{<1>}}\xi_{H}^{-2}(m_{1}))S^{-1}(b_{[-1]})\>
  \<f,(a\cdot(\xi_{A}(b_{[0]_{<0>}})\cdot\xi_{M}^{-1}(m_{0})))_{1}\>\\
  &&\xi_{M}^{2}((\xi_{A}(a)\cdot\xi_{M}(\xi_{A}(b_{[0]_{<0>}})\cdot\xi_{M}^{-1}(m_{0})))_{0})\\
  &\stackrel{(\ref{eq10})}{=}&
  \<g,(b\cdot\xi_{M}^{-1}(m))_{1}\>
  \<\xi_{H}^{*-1}(f),(\xi_{A}(a)\cdot\xi_{M}((b\cdot\xi_{M}^{-1}(m))_{0}))_{1}\>\\
  && \xi_{M}^{2}((\xi_{A}(a)\cdot\xi_{M}((b\cdot\xi_{M}^{-1}(m))_{0}))_{0})\\
  &=&(\xi_{H}^{*-1}(f)\bowtie \xi_{A}(a))\triangleright\xi_{M}^{2}((b\cdot\xi_{M}^{-1}(m))_{0})
  \<g,(b\cdot\xi_{M}^{-1}(m))_{1}\>\\
  &=&(\xi_{H}^{*-1}(f)\bowtie \xi_{A}(a))\triangleright[(g\bowtie b)\triangleright m],
  \end{eqnarray*}
  as needed. It is not hard to see that
  $(\varepsilon_{H}\bowtie 1_{A})\triangleright m=\xi_{M}(m)$,
   for all $m\in M,$ so $M$ is a left $H^{*}\bowtie A$-module.
   The fact that a morphism in $ _{A}\mathcal{MHYD}^{H}(H)$
   becomes a morphism in $_{H^{*}\bowtie A}\mathcal{M}$
   can be proved more easily, we leave the details to the reader.
   \hfill $\blacksquare$
   \\

{\bf Lemma 2.6.} Let $(H,\xi_{H})$ be a monoidal Hom-Hopf algebra
and $(H,A,H)$ a Yetter-Drinfeld Hom-datum and assume $H$ is finite dimensional.
 We have a functor $G:  _{H^{*}\bowtie A}\mathcal{M}\rightarrow_{A}\mathcal{MHYD}^{H}(H),$
given by $G(M)=M$ as $k$-module, with the structure maps
defined by
 \begin{eqnarray}\label{2.7}
  u\cdot m=(\varepsilon_{H}\bowtie\xi_{A}^{-1}(u))\triangleright m,
  \end{eqnarray}
  \begin{eqnarray}\label{2.8}
  \rho_{M}:M\rightarrow M\o H,\ \ \rho_{M}(m)= m_{0}\o m_{1}
  =\sum_{i=1}^{n}(\xi_{H}^{*2}(e^{i})\bowtie 1_{A})\triangleright \xi_{M}^{-2}(m)\o e_{i},
  \end{eqnarray}
  for all $u\in (A,\xi_{A})$ and $m\in (M,\xi_{M})$. Here $\{e_{i}\}_{i=1,...,n}$
  is a basis of $H$ and $\{e^{i}\}_{i=1,...,n}$ is the corresponding dual basis
  of $H^{*}$. $G$ transforms a morphism to itself.
  \\

{\bf Proof.} The most difficult part of the proof is to show that $G(M)$
 satisfies the relations (\ref{eq10}) or (\ref{eq11}). It is then straightforward to show that
 a map in $ _{H^{*}\bowtie A}\mathcal{M}$ is also a map in $_{A}\mathcal{MHYD}^{H}(H),$
 and that $G$ is a functor.

 We compute:
  \begin{eqnarray*}
  &&u_{<0>}\cdot m_{0}\o u_{<1>}m_{1}\\
  &=& \sum_{i=1}^{n}(\varepsilon_{H}\bowtie \xi_{A}^{-1}(u_{<0>}))
  \triangleright((\xi_{H}^{*2}(e^{i})\bowtie 1_{A})\triangleright
  \xi_{M}^{-2}(m))\o u_{<1>}e_{i}\\
  &=&\sum_{i=1}^{n}[\varepsilon_{H}(\xi_{A}^{-2}(u_{<0>_{[-1]}})
  \rightharpoonup(\xi_{H}^{*4}(e^{i})\leftharpoonup S^{-1}\xi_{H}^{-2}(u_{<0>_{[0]_{<1>}}})))
  \bowtie \xi_{A}(u_{<0>_{[0]_{<0>}}})]\triangleright\xi_{M}^{-1}(m)\\
  &&\o u_{<1>}e_{i}\\
  &=&\sum_{i=1}^{n}\<e_{1}^{i},S^{-1}\xi_{H}(u_{<0>_{[0]_{<1>}}})\>\<e_{22}^{i},\xi_{H}^{-1}(u_{<0>_{[-1]}})\>
  \<e_{21}^{i},\xi_{H}^{-1}((\xi_{A}(u_{<0>_{[0]_{<0>}}})\cdot\xi_{M}^{-2}(m))_{1})\>\\
  &&\xi_{M}^{2}((\xi_{A}(u_{<0>_{[0]_{<0>}}})\cdot\xi_{M}^{-2}(m))_{0})
  \o u_{<1>}e_{i}\\
  &=&\sum_{i=1}^{n}\<e^{i},S^{-1}\xi_{H}(u_{<0>_{[0]_{<1>}}})((u_{<0>_{[0]_{<0>}}}\cdot\xi_{M}^{-3}(m))_{1}
  \xi_{H}^{-1}(u_{<0>_{[-1]}}))\>\\
  &&\xi_{M}^{2}((\xi_{A}(u_{<0>_{[0]_{<0>}}})\cdot\xi_{M}^{-2}(m))_{0})
  \o u_{<1>}e_{i}\\
  &=&\xi_{M}^{2}((\xi_{A}(u_{<0>_{[0]_{<0>}}})\cdot\xi_{M}^{-2}(m))_{0})
  \o (\xi_{H}^{-1}(u_{<1>})S^{-1}\xi_{H}(u_{<0>_{[0]_{<1>}}}))\\
  &&  \xi_{H}((u_{<0>_{[0]_{<0>}}}\cdot\xi_{M}^{-3}(m))_{1}u_{<0>_{[-1]}})\\
  &=&\xi_{M}^{2}((u_{<0>_{[0]}}\cdot\xi_{M}^{-2}(m))_{0})
  \o (u_{<1>2})S^{-1}(u_{<1>1}))
   ((u_{<0>_{[0]}}\cdot\xi_{M}^{-2}(m))_{1}\xi_{H}(u_{<0>_{[-1]}}))\\
  &=&\xi_{M}((u_{[0]}\cdot\xi_{M}^{-1}(m))_{0})\o (u_{[0]}\cdot\xi_{M}^{-1}(m))_{1}u_{[-1]},
  \end{eqnarray*}
  for all $u\in (A,\xi_{A})$ and $m\in(M,\xi_{M})$, and this finishes the proof.
   \hfill $\blacksquare$\\

   The next result generalizes (\cite{CZ2014}, Prop. 4.3),
    which is recovered by taking $H = A$.
    \\

   {\bf Theorem 2.7.} Let $(H,\xi_{H})$ be a monoidal Hom-Hopf algebra
and $(H,A,H)$ a Yetter-Drinfeld datum, assuming $H$ to be finite dimensional.
 Then the categories $_{A}\mathcal{MHYD}^{H}(H)$ and $ _{H^{*}\bowtie A}\mathcal{M}$
are isomorphic.  \\

{\bf Proof.} We have to verify that the functors $F$ and $G$ defined in Lemmas 2.5 and 2.6
are inverse to each other. Let $M\in _{A}\mathcal{MHYD}^{H}(H)$. The structures on
$G(F(M))$ are denoted by $\cdot^{'}$
and $\rho_{M}^{'}$. For any $u\in (A,\xi_{A})$ and $m\in (M,\xi_{M})$ we have that
$$u\cdot^{'}m = (\v\bowtie \xi_{A}^{-1}(u))\triangleright m
= \<\v,(\xi_{A}^{-1}(u)\cdot \xi_{M}^{-1}(m))_{1}\>
\xi_{M}^{2}((\xi_{A}^{-1}(u)\cdot \xi_{M}^{-1}(m))_{0}) = u\cdot m.$$
We now compute for $m\in (M,\xi_{M})$ that
 \begin{eqnarray*}
  \rho_{M}^{'}(m)&=&\sum_{i=1}^{n}(\xi_{H}^{*2}(e^{i})\bowtie1_{A})
  \triangleright \xi_{M}^{-2}(m)\o e_{i}\\
  &\stackrel{(\ref{2.6})}{=}&\sum_{i=1}^{n}\<\xi_{H}^{*2}(e^{i}),(1_{A}\cdot\xi_{M}^{-3}(m))_{1}\>
  \xi_{M}^{2}((1_{A}\cdot\xi_{M}^{-3}(m))_{0})\o e_{i}\\
  &=&\sum_{i=1}^{n}\<e^{i},m_{1}\>m_{0}\o e_{i} = \rho_{M}(m).
  \end{eqnarray*}

  Conversely, take $M\in  _{H^{*}\bowtie A}\mathcal{M}$. We want to show that
  $F(G(M))=M$. If we denote the left $H^{*}\bowtie A$-action on $F(G(M))$ by
  $\mapsto$, then using Lemmas 2.5 and 2.6 we find,
  for all $f\in (H^{*},\xi_{H}^{*-1}),u\in (A,\xi_{A})$ and $m\in (M,\xi_{M})$:
  \begin{eqnarray*}
  (f\bowtie u)\mapsto m &=&\<f,(u\cdot\xi_{M}^{-1}(m))_{1}\>\xi_{M}^{2}((u\cdot\xi_{M}^{-1}(m))_{0})\\
  &=&\sum_{i=1}^{n}\<f,e_{i}\>\xi_{M}^{2}((\xi_{H}^{*2}(e^{i})\bowtie 1_{A})
  \triangleright \xi_{M}^{-2}(u\cdot\xi_{M}^{-1}(m)))\\
  &=&\<\xi_{H}^{*2}(f),\xi_{M}^{-2}(u\cdot\xi_{M}^{-1}(m))\>\xi_{M}^{2}(u\cdot\xi_{M}^{-1}(m))\\
  &=&(f\bowtie u)\triangleright m,
  \end{eqnarray*}
and this finishes our proof.  \hfill $\blacksquare$\\

{\bf Proposition 2.8.} Let  $(H,\xi_{H})$ be finite dimensional
and $H(\a,\b)$ be an $H$-Hom-bicomdule algebra, with an
$H$-Hom-comodule structures showed in Example 2.9 (in \cite{YW2014}).
Then $_{H(\a,\b)}\mathcal{MHYD}^{H}(H)\simeq  _{H^{*}\bowtie H(\a,\b)}\mathcal{M}$.
 \\

The proof is left to the reader.
\\

 Recall from Prop.2.12 in \cite{YW2014},
$_{H}\mathcal{MHYD}^{H}(\a,\b)=_{H(\a,\b)}\mathcal{MHYD}^{H}(H)$.
\\

 {\bf Proposition 2.9.} $_{H}\mathcal{MHYD}^{H}(\a,\b)
 \simeq  _{H^{*}\bowtie H(\a,\b)}\mathcal{M}$.
 \\

 We just give the correspondence as follows. If $M\in _{H}\mathcal{MHYD}^{H}(\a,\b)$,
then $M\in  _{H^{*}\bowtie H(\a,\b)}\mathcal{M}$ with structure
$$ (f\bowtie h)\triangleright m=f((h\cdot\xi_{M}^{-1}(m))_{1})
  \xi_{M}^{2}((h\cdot\xi_{M}^{-1}(m))_{0}).$$
  Conversely, if $M\in  _{H^{*}\bowtie H(\a,\b)}\mathcal{M}$, then
$M\in _{H}\mathcal{MHYD}^{H}(\a,\b)$ with structures
\begin{eqnarray*}
   h\cdot m &=& (\varepsilon_{H}\bowtie\xi_{H}^{-1}(h))\triangleright m,\\
 \rho_{M}(m)&=& m_{0}\o m_{1}
  =(\sum_{i=1}^{n}\xi_{H}^{*2}(e^{i})\bowtie 1_{A})\triangleright \xi_{M}^{-2}(m)\o e_{i}
  \end{eqnarray*}
 for all $f\in H^{*},h\in H,m\in M,$
 where $\{e_{i}\}_{1,...,n},\{e^{i}\}_{1,...,n}$ are dual bases in $H$ and $H^{*}$.
The proof is left to the reader.

\section*{3. A BRAIDED $T$-CATEGORY $Rep(\mathcal{MHD}(H))$ }
\def\theequation{3. \arabic{equation}}
\setcounter{equation} {0} \hskip\parindent

  Denote $G={\sl Aut}_{mHH}(H) \times {\sl Aut}_{mHH}(H)$
 a group with multiplication as follows:
  for all $\a,\b,\g,\d\in {\sl Aut}_{mHH}(H)$,
$$
(\a,\b)\ast (\g,\d)=(\a\g, \d\g^{-1}\b\g).
$$
 The unit of this group is $(id,id)$ and $(\a,\b)^{-1}=(\a^{-1}, \a\b^{-1}\a^{-1})$.
\\

 In this section we will construct a monoidal Hom-Hopf $T$-coalgebra
 over $G$, denoted by $\mathcal{MHD}(H)$, and prove that the $T$-category
 $Rep(\mathcal{MHD}(H))$ of representation of $\mathcal{MHD}(H)$
 is isomorphic to $\mathcal {MHYD}(H)$ in \cite{YW2014} as braided $T$-categories.
\\


 {\bf Proposition 3.1.} Let $(M,\xi_{M})\in _{H}\mathcal{MHYD}^{H}(\a,\b)$
 and assume that $(M,\xi_{M})$ is finite dimensional.
 Then $(M^{*},\xi_{M}^{*-1})=Hom(M,k)$ becomes an object in
  $_{H}\mathcal{MHYD}^{H}(\a^{-1},\a\b^{-1}\a^{-1})$, with module structure
  $$ (h\cdot p)(m)=p(\b^{-1}\a^{-1}S\xi_{H}^{-1}(h)\cdot\xi_{M}^{-2}(m)),$$
  and comodule structure
  $$\rho(p)(m)=  p_{0}(\xi_{M}^{-1}(m))\o\xi_{H}(p_{1})
  = p(\xi_{M}(m_{0}))\o S^{-1}\xi_{H}^{2}(m_{1}),$$
  for all $h\in H,p\in M^{*}$ and $m\in M$. Moreover, the maps
  $b_{M}:k\rightarrow M\o M^{*}, ~b_{M}(1)=\sum_{i}c_{i}\o c^{i}$
  (where $\{c_{i}\}$ and $\{c^{i}\}$ are dual bases in $M$ and $M^{*}$)
  and $d_{M}:M^{*}\o M\rightarrow k, ~d_{M}(p\o m)=p(m),$ are left
  $H$-module maps and right $H$-comodule maps and we have
  $$(\xi_{M}\o d_{M})(b_{M}\o \xi_{M}^{-1})=id_{M},\ \
  ( d_{M}\o \xi_{M}^{*-1})(\xi_{M}^{*}\o b_{M})=id_{M^{*}}.$$
\\

  {\bf Proof.} Following the idea of the proof of Panaite and Staic
  (\cite{PS2007}, Prop. 3.6), we first prove that $(M^{*},\xi_{M}^{*-1})$
  is indeed an object in $_{H}\mathcal{MHYD}^{H}(\a^{-1},\a\b^{-1}\a^{-1})$.
  We compute:
 \begin{eqnarray*}
  && (\xi_{H}(h_{21})\cdot p_{0})(m)\o (\a\b^{-1}\a^{-1}(h_{22})
  \xi_{H}^{-1}(p_{1})) \a^{-1}S^{-1}(h_{1})\\
  &=& p_{0}(\b^{-1}\a^{-1}S(h_{21})\cdot\xi_{M}^{-2}(m)
  \o (\a\b^{-1}\a^{-1}(h_{22})\xi_{H}^{-1}(p_{1}))\a^{-1}S^{-1}(h_{1})\\
  &=& p(\xi_{M}^{2}((\b^{-1}\a^{-1}S(h_{21})\cdot\xi_{M}^{-2}(m))_{0}))
  \o (\a\b^{-1}\a^{-1}(h_{22})S^{-1}((\b^{-1}\a^{-1}S(h_{21})\cdot\xi_{M}^{-2}(m))_{1}))\\
  && \a^{-1}S^{-1}(h_{1})\\
  &=& p(\xi_{M}^{2}(\b^{-1}\a^{-1}S\xi_{H}(h_{2112})\cdot\xi_{M}^{-2}(m_{0})))
  \o (\a\b^{-1}\a^{-1}(h_{22})S^{-1}((\a^{-1}S(h_{2111})\xi_{H}^{-3}(m_{1}))\\
  && \a\b^{-1}\a^{-1}(h_{212})))\a^{-1}S^{-1}(h_{1})\\
  &=& p(\b^{-1}\a^{-1}S\xi_{H}^{3}(h_{2112})\cdot m_{0})
  \o (\a\b^{-1}\a^{-1}(h_{22})(\a\b^{-1}\a^{-1}S^{-1}(h_{212})\\
  && (S^{-1}\xi_{H}^{-3}(m_{1})\a^{-1}(h_{2111}))))\a^{-1}S^{-1}(h_{1})\\
  &=& p(\b^{-1}\a^{-1}S\xi_{H}^{3}(h_{2112})\cdot m_{0})
  \o ((\a\b^{-1}\a^{-1}\xi_{H}^{-1}(h_{22})\a\b^{-1}\a^{-1}S^{-1}(h_{212}))\\
  && (S^{-1}\xi_{H}^{-2}(m_{1})\a^{-1}\xi_{H}(h_{2111})))\a^{-1}S^{-1}(h_{1})\\
   &=& p(\b^{-1}\a^{-1}S\xi_{H}^{3}(h_{2112})\cdot m_{0})
  \o (\a\b^{-1}\a^{-1}(h_{22})\a\b^{-1}\a^{-1}S^{-1}\xi_{H}(h_{212}))\\
  && ((S^{-1}\xi_{H}^{-2}(m_{1})\a^{-1}\xi_{H}(h_{2111}))\a^{-1}S^{-1}\xi_{H}^{-1}(h_{1}))\\
   &=& p(\b^{-1}\a^{-1}S\xi_{H}(h_{21})\cdot m_{0})
  \o (\a\b^{-1}\a^{-1}\xi_{H}(h_{222})\a\b^{-1}\a^{-1}S^{-1}\xi_{H}(h_{221}))\\
  && (S^{-1}\xi_{H}^{-1}(m_{1})(\a^{-1}\xi_{H}^{-1}(h_{12})\a^{-1}S^{-1}\xi_{H}^{-1}(h_{11})))\\
  &=& p(\b^{-1}\a^{-1}S\xi_{H}^{-1}(h)\cdot m_{0}) \o S^{-1}\xi_{H}(m_{1})\\
  &=& (h\cdot p)(\xi_{M}^{2}(m_{0}))\o S^{-1}\xi_{H}(m_{1})\\
  &=& (h\cdot p)_{0}(m)\o (h\cdot p)_{1},
  \end{eqnarray*}
  which means that
  \begin{eqnarray*}
  (h\cdot p)_{0}\o (h\cdot p)_{1}
  = (\xi_{H}(h_{21})\cdot p_{0})\o (\a\b^{-1}\a^{-1}(h_{22})
  \xi_{H}^{-1}(p_{1})) \a^{-1}S^{-1}(h_{1}).
  \end{eqnarray*}

  On $k$ we have the trivial Hom-module and Hom-comodule structure,
  and with these $k\in _{H}\mathcal{YD}^{H}$. We want to prove that $b_{M}$
  and $d_{M}$ are $H$-Hom-module maps. We compute:
  \begin{eqnarray*}
   (h\cdot b_{M}(1))(m)&=&(h\cdot(\sum_{i}c_{i}\o c^{i}))(m)\\
   &\stackrel{(\ref{eq12})}{=}& \sum_{i}\a^{-1}(h_{1})\cdot c_{i}\o (\a\b\a^{-1}(h_{2})
   \cdot c^{i})(m)\\
   &=& \sum_{i}\a^{-1}(h_{1})\cdot c_{i}\o c^{i}(\b^{-1}\a^{-1}S\a\b\a^{-1}\xi_{H}^{-1}(h_{2})
   \cdot \xi_{M}^{-2}(m))\\
    &=& \sum_{i}\a^{-1}(h_{1})\cdot c_{i}\o c^{i}(S\a^{-1}\xi_{H}^{-1}(h_{2})\cdot \xi_{M}^{-2}(m))\\
    &=& \a^{-1}(h_{1})\cdot (S\a^{-1}\xi_{H}^{-1}(h_{2})\cdot \xi_{M}^{-2}(m))\\
    &=& \a^{-1}(\xi_{-1}(h_{1})S\xi_{H}^{-1}(h_{2}))\cdot \xi_{M}^{-1}(m)\\
    &=& \varepsilon(h)\sum_{i}c_{i}\o c^{i}(m)\\
    &=& (\varepsilon(h)b_{M}(1))(m),
  \end{eqnarray*}
  \begin{eqnarray*}
   d_{M}(h\cdot(p\o m))&\stackrel{(\ref{eq12})}{=}&
   d_{M}(\a(h_{1})\cdot p\o \b^{-1}(h_{2})\cdot m)\\
   &=& (\a(h_{1})\cdot p)(\b^{-1}(h_{2})\cdot m)\\
   &=& p(\b^{-1}\a^{-1}S\a\xi_{H}^{-1}(h_{1})\cdot\xi_{M}^{-2}(\b^{-1}(h_{2})\cdot m))\\
   &=& p(\b^{-1}(S\xi_{H}^{-2}(h_{1})\xi_{H}^{-2}(h_{2}))\cdot\xi_{M}^{-1}( m))\\
   &=& \varepsilon(h)d_{M}(p\o m).
  \end{eqnarray*}
  They also are $H$-Hom-comodule maps;
  \begin{eqnarray*}
   ((b_{M}(1))_{0}\o (b_{M}(1))_{1})(m)
   &=& \sum_{i}(c_{i})_{0}\o (c^{i})_{0}(m)\o(c^{i})_{1}(c_{i})_{1}\\
   &=& \sum_{i}(c_{i})_{0}\o c^{i}(\xi_{M}^{2}(m_{0}))\o S^{-1}\xi_{H}(m_{1})(c_{i})_{1}\\
   &=& \xi_{M}^{2}(m_{00})\o S^{-1}\xi_{H}(m_{1})\xi_{H}^{2}(m_{01})\\
   &=& \xi_{M}(m_{0})\o S^{-1}\xi_{H}^{2}(m_{12})\xi_{H}^{2}(m_{11})\\
   &=& (b_{M}(1)\o 1)(m),
  \end{eqnarray*}
  \begin{eqnarray*}
   d_{M}((p\o m)_{0})\o (p\o m)_{1}&=& p_{0}( m_{0})\o m_{1}p_{1}\\
   &=& p(\xi_{M}^{2}(m_{00}))\o m_{1}S^{-1}\xi_{H}(m_{01})\\
   &=& p(\xi_{M}(m_{0}))\o \xi_{H}(m_{12})S^{-1}\xi_{H}(m_{11})\\
   &=& d_{M}(p\o m)\o 1.
  \end{eqnarray*}
  Finally, we compute:
   \begin{eqnarray*}
   (\xi_{M}\o d_{M})(b_{M}\o \xi_{M}^{-1})(m)
   &=& (\xi_{M}\o d_{M})(b_{M}(1)\o \xi_{M}^{-1}(m))\\
   &=& (\xi_{M}\o d_{M})(\sum_{i}(c_{i}\o c^{i})\o \xi_{M}^{-1}(m))\\
   &=& \sum_{i}\xi_{M}^{2}(c_{i})\o c^{i}(\xi_{M}^{-2}(m)) = m
  \end{eqnarray*}
  The argument for $( d_{M}\o \xi_{M}^{*-1})(\xi_{M}^{*}\o b_{M})=id_{M^{*}}$
  is analogous.  \hfill $\blacksquare$\\

  Similarly, one can obtain:
  \\

  {\bf Proposition 3.2.} Let $(M,\xi_{M})\in _{H}\mathcal{MHYD}^{H}(\a,\b)$
 and assume that $(M,\xi_{M})$ is finite dimensional.
 Then $(^{*}M,^{*}\xi_{M}^{-1})=Hom(M,k)$ becomes an object in
  $_{H}\mathcal{MHYD}^{H}(\a^{-1},\a\b^{-1}\a^{-1})$, with module structure
  $$ (h\cdot p)(m)=p(\b^{-1}\a^{-1}S\xi_{H}^{-1}(h)\cdot\xi_{M}^{-2}(m)),$$
  and comodule structure
  $$\rho(p)(m)=  p_{0}(\xi_{M}^{-1}(m))\o\xi_{H}(p_{1})
  = p(\xi_{M}(m_{0}))\o S^{-1}\xi_{H}^{2}(m_{1}),$$
  for all $h\in H, ~p\in\ ^{*}M$ and $m\in M$. Moreover, the maps
  $b_{M}:k\rightarrow M\o\ ^{*}M, ~b_{M}(1)=\sum_{i}c_{i}\o c^{i}$
  (where $\{c_{i}\}$ and $\{c^{i}\}$ are dual bases in $M$ and $^{*}M$)
  and $d_{M}:\ ^{*}M\o M\rightarrow k, ~d_{M}(p\o m)=p(m),$ are left
  $H$-module maps and right $H$-comodule maps and we have
  $$(\xi_{M}\o d_{M})(b_{M}\o \xi_{M}^{-1})=id_{M},\ \
  ( d_{M}\o ^{*}\xi_{M}^{-1})(^{*}\xi_{M}\o b_{M})=id_{^{*}M}.$$
\\

  Now, if we consider $\mathcal {MHYD}(H)_{fd}$, the subcategory
  of $\mathcal {MHYD}(H)$ consisting of finite dimensional objects,
  then by Proposition 3.1. and Proposition 3.2. we obtain:\\

  {\bf Corollary 3.3.} $\mathcal {MHYD}(H)_{fd}$ is a braided $T$-category
  with left and right dualities over $G$.
  \\

  Assume now that $(H,\xi_{H})$ is finite dimensional. We will
  construct a monoidal Hom-Hopf $T$-coalgebra over $G$, denoted by
  $ \mathcal{MHD}(H)$, with the property that the $T$-category
  $Rep(\mathcal{MHD}(H))$ of representation of $\mathcal{MHD}(H)$ is isomorphic
   to $\mathcal{MHYD}(H)$ as braided $T$-categories.
  \\

  {\bf Theorem 3.4.} $\mathcal{MHD}(H)=\{\mathcal{MHD}(H)_{(\a,\b)}\}_{(\a,\b)\in G}$
  is a monoidal Hom-Hopf $T$-coalgebra with the following structures:
  \begin{itemize}
   \item For any $(\a,\b)\in G,$ the $(\a,\b)$-component $\mathcal{MHD}(H)_{(\a,\b)}$
   will be the diagonal crossed Hom-product algebra $H^{*}\bowtie H(\a,\b)$
    in Eq. (\ref{2.10}),
   \item The comultiplication on $ \mathcal{MHD}(H)$ is given by
   \begin{eqnarray*}
   \Delta_{(\a,\b),(\g,\d)}:\mathcal{MHD}(H)_{(\a,\b)*(\g,\d)}
   &\rightarrow& \mathcal{MHD}(H)_{(\a,\b)}\o\mathcal{MHD}(H)_{(\g,\d)},\\
   \Delta_{(\a,\b),(\g,\d)}(f\bowtie h)&= &(f_{1}\bowtie\g(h_{1}))\o(f_{2}\bowtie\g^{-1}\b\g(h_{2})),
   \end{eqnarray*}
   \item The counit $\varepsilon$ is obtained by setting
   $$\varepsilon(f\bowtie h)=\varepsilon(h)f(1_{H}),$$
   \item For any $(\a,\b)\in G$, the $(\a,\b)^{th}$ component of the antipode of
   $\mathcal{MHD}(H)$ is given by
   \begin{eqnarray*}
   S_{(\a,\b)}:\mathcal{MHD}(H)_{(\a,\b)}&\rightarrow& \mathcal{MHD}(H)_{(\a,\b)^{-1}}= \mathcal{MHD}(H)_{(\a^{-1},\a\b^{-1}\a^{-1})},\\
   S_{(\a,\b)}(f\bowtie h)&=&(\varepsilon\bowtie\a\b S\xi_{H}^{-1}(h))
   (S^{*-1}\xi_{H}^{*}(f)\bowtie1_{H}),
   \end{eqnarray*}
    \item For $(\a,\b),(\g,\d)\in G$, the conjugation isomorphism is given by
    \begin{eqnarray*}
    \varphi_{(\g,\d)}^{(\a,\b)}: \mathcal{MHD}(H)_{(\g,\d)}
    &\rightarrow& \mathcal{MHD}(H)_{(\a,\b)*(\g,\d)*(\a,\b)^{-1}},\\
    \varphi_{(\g,\d)}^{(\a,\b)}(f\bowtie h)&=&(f\circ\b\a^{-1}
    \bowtie\a\g^{-1}\b^{-1}\g(h)),
   \end{eqnarray*}
  \end{itemize}
for all $f\in H^{*},h\in H.$
\\

 {\bf Proof.} We have to check the axioms of monoidal Hom-Hopf $T$-coalgebra.
 Hom-coassociativity and multiplicativity of $\Delta$ are satisfied.
 We compute
  \begin{eqnarray*}
    && m_{(\a,\b)}(id\o S_{(\a,\b)^{-1}})\Delta_{(\a,\b),(\a,\b)^{-1}}(f\bowtie h)\\
    &=& (f_{1}\bowtie\a^{-1}(h_{1}))[(\varepsilon\bowtie S\a^{-1}\xi_{H}^{-1}(h_{2}))
    (S^{*-1}\xi^{*}(f_{2})\bowtie 1_{H})]\\
    &=& [(\xi_{H}^{*}(f_{1})\bowtie\a^{-1}\xi_{H}^{-1}(h_{1}))
    (\varepsilon\bowtie S\a^{-1}\xi_{H}^{-1}(h_{2}))]
    (S^{*-1}(f_{2})\bowtie 1_{H})\\
    &=& [(\xi_{H}^{*}(f_{1})(\a((\a^{-1}\xi_{H}^{-1}(h_{1}))_{1})\rightharpoonup
    (\varepsilon\leftharpoonup S^{-1}\b((\a^{-1}\xi_{H}^{-1}(h_{1}))_{22})))\\
    && \bowtie\xi_{H}^{2}((\a^{-1}\xi_{H}^{-1}(h_{1}))_{21})
     S\a^{-1}\xi_{H}^{-1}(h_{2}))](S^{*-1}(f_{2})\bowtie 1_{H})\\
     &=& [f_{1} \varepsilon( S^{-1}\b\a^{-1}\xi_{H}^{-2}(h_{122}))
     \varepsilon(\xi_{H}^{-2}(h_{11}))\bowtie\a^{-1}\xi_{H}(h_{121})
     S\a^{-1}\xi_{H}^{-1}(h_{2})]\\
     &&(S^{*-1}(f_{2})\bowtie 1_{H})\\
     &=& (f_{1} \bowtie\a^{-1}\xi_{H}^{-1}(h_{1}S(h_{2})))(S^{*-1}(f_{2})\bowtie 1_{H})\\
     &=& \varepsilon(h)(f_{1} \bowtie 1_{H})(S^{*-1}(f_{2})\bowtie 1_{H})\\
     &=& \varepsilon(h)(f_{1}S^{*-1}(f_{2})\bowtie 1_{H})\\
     &=& \varepsilon(f\bowtie h)(\varepsilon\bowtie 1_{H}),\\
   \end{eqnarray*}
   and similarly
   $ m_{(\a,\b)}( S_{(\a,\b)^{-1}}\o id)\Delta_{(\a,\b)^{-1},(\a,\b)}(f\bowtie h)
   = \varepsilon(f\bowtie h)(\varepsilon\bowtie 1_{H}).$
   This completes the proof.  \hfill $\blacksquare$\\

 Moreover, via the isomorphisms Prop.2.9. $_{H}\mathcal{MHYD}^{H}(\a,\b)
 \simeq  _{H^{*}\bowtie H(\a,\b)}\mathcal{M}$, we obtain
  \\

  {\bf Theorem 3.5.} $Rep(\mathcal{MHD}(H))\simeq \mathcal{MHYD}(H)$
  as braided $T$-categories over $G$.

\section*{4. A BRAIDED $T$-CATEGORY $\mathcal{ZMHYD}(H)$}
\def\theequation{4. \arabic{equation}}
\setcounter{equation} {0} \hskip\parindent

In this section, we will construct a new braided $T$-category
$\mathcal{ZMHYD}(H)$ over $\mathbb{Z}.$
\\

 {\bf Definition 4.1.} Let $(C,\xi_{C})$ be a monoidal Hom-coalgebra.
 Then $g$ is called a group-like element, that is
 $$\xi_{C}(g) = g,\ \ \Delta(g)= g \o g,\ \  \varepsilon(g)=1, $$
 for all $g\in C$.
\\

 {\bf Example 4.2.} Recall from Example 3.5 in \cite{CWZ2013}
 that $( H_{4} =k\{ 1,\ g,\ x, \ y = gx\,\},\xi_{H_{4}} , \Delta , \varepsilon , S )$
 is a monoidal Hom-Hopf algebra, where the algebraic structure are given as follows:

$\bullet$   The multiplication
 $"\circ "$ is given by

 $$\begin{array}{|c|c|c|c|c|}
  \hline  \circ & 1_{H_4} & g & x & y\\
 \hline 1_{H_4} & 1_{H_4} &  g & cx & cy \\
  g & g & 1_{H_4} & cy & cx \\
  x& cx & -cy& 0 & 0 \\
  y& cy & -cx & 0 & 0 \\
 \hline
 \end{array}$$

$\bullet$ The automorphism $\xi_{H_{4}}$ is given by
 $ \xi_{H_{4}}(1)=1,\ \xi_{H_{4}}(g)=g,\ \xi_{H_{4}}(x)=cx,$\\
  $\xi_{H_{4}}(gx)=cgx,$
  for all $0\neq c\in k;$

$\bullet$  The comultiplication $\Delta$ is defined by
 \begin{eqnarray*}
 \Delta(1)=1\otimes 1,&&
  \Delta(g)=g\otimes g,\\
 \Delta(x)=c^{-1}(x\otimes 1)+ c^{-1}(g\otimes x),
 &&\Delta(gx)=c^{-1}(gx\otimes g)+c^{-1}(1\otimes gx);
 \end{eqnarray*}

$\bullet$ The counit $\varepsilon$ is defined by
 $ \varepsilon(1)=1,\ \
 \varepsilon(g)=1,\ \
 \varepsilon(x)=0,\ \  \varepsilon(gx)=0. $

$\bullet$ The antipode $S$ is given by
 $ S(1)=1,\ \  S(g)=g,\ \   S(x)=-gx,\ \  S(gx)=-x.$

  Then $1_{H_{4}}, g $ are group-like elements of $H_{4}$.
 \\

 In \cite{YW2014}, the authors introduced the notion of
 left-right $(\a, \b)$-Yetter-Drinfeld Hom-module {\it Definition 1.2.6.}
 We will in the section give its some special cases.
 \\

 {\bf Example 4.3.} For $(M,\xi_{M})\in _{H}\mathcal{MHYD}^{H}(S^{2},id)$,
   the {\it left-right anti-Yetter-Drinfeld Hom-module category, }
    i.e., the compatibility condition is
     \begin{eqnarray*}
  (h\c m)_{0}\o (h\c m)_{1} =
  \xi_{H}(h_{21})\c m_{0}\o (h_{22}\xi_{H}^{-1}(m_{1}))S(h_{1}),
  \end{eqnarray*}
  for $h\in H,m\in M.$

   {\bf Example 4.4.} In $_{H}\mathcal{MHYD}^{H}(S^{2n},id)$, the object is called
   a {\it left-right $n$-Yetter-Drinfeld Hom-modules, i.e., $n-\mathcal{MHYD}$-module},
   for $(M,\xi_{M})\in _{H}\mathcal{MHYD}^{H}(S^{2n},id)$,
    the compatibility condition is
     \begin{eqnarray*}
  (h\c m)_{0}\o (h\c m)_{1} =
  \xi_{H}(h_{21})\c m_{0}\o (h_{22}\xi_{H}^{-1}(m_{1}))S^{2n-1}(h_{1}),
  \end{eqnarray*}
  for $h\in H,m\in M.$

{\bf Example 4.5.} Similar to Panaite and Staic (\cite{PS2007}, Example 2.7),
for $\a,\b\in {\sl Aut}_{mHH}(H) $, and assume that there is an algebra map
$\theta: H\rightarrow k$ and a group-like element $\omega\in (H,\xi_{H})$ such that
\begin{eqnarray*}
  \a(h) = \om^{-1}(\theta(h_{11})\b(h_{12})\theta(S(h_{2}))\om), \ \ \ \forall h\in H.
  \end{eqnarray*}

  Then we can check that $k\in  _{H}\mathcal{MHYD}^{H}(\a,\b)$
  with structures: $h\cdot 1=\theta(h)$ and $\r(1)= 1\o \om.$
  More generally, if $V$ is any vector space,
  then $(V,\xi_{V})\in _{H}\mathcal{MHYD}^{H}(\a,\b),$
  with structures $h\cdot v=\th(h)\xi_{V}(v)$ and
  $\r(v)= v_{0}\o v_{1}= \xi_{V}^{-1}(v)\o \om,$ for all
  $h\in H$ and $v\in V.$
  \\

  If $\a,\b\in {\sl Aut}_{mHH}(H)$ such that there exist $(\th,\om)$
  as show in Example 4.1, we will say that $(\th,\om)$ is a pair
  in involution corresponding to $(\a,\b)$ and the left-right
  $(\a,\b)$-Yetter-Drinfeld Hom-modules $k$ and $(V,\xi_{V})$ will be denoted
  by $_{\th}k^{\om}$ and $_{\th}V^{\om}$, respectively.
  \\


  In the following, we will show that in the presence of a pair
  in involution, there exists an isomorphism of categories
  $_{H}\mathcal{MHYD}^{H}(\a,\b)\simeq _{H}\mathcal{MHYD}^{H}.$
\\

{\bf Proposition 4.6.} Let $\a,\b\in {\sl Aut}_{mHH}(H)$ and assume
that there exists $(\th,\om)$ a pair in involution corresponding to $(\a,\b)$.
Then the categories $_{H}\mathcal{MHYD}^{H}(\a,\b)$
and $_{H}\mathcal{MHYD}^{H}$ are isomorphic.
\\

 {\bf Proof.} In order to prove the isomorphism between
 two categories, we only need to give a pair of inverse functors.
 The functors pair $(F,G)$ is given as follows.

 If $(M,\xi_{M})\in _{H}\mathcal{MHYD}^{H}(\a,\b),$
  then $(F(M),\xi_{M})\in _{H}\mathcal{MHYD}^{H}$, where $F(M)= M$
  as vector space, with structures
  \begin{eqnarray*}
  &h\rightarrow m = &\th(\b^{-1}S(h_{1}))\b^{-1}(h_{2})\cdot m, \\
  &\r(m)=:& m_{<0>}\o m_{<1>}=m_{0}\o m_{1}\om^{-1}.
  \end{eqnarray*}

   If $(N,\xi_{N})\in _{H}\mathcal{MHYD}^{H},$
  then $(G(N),\xi_{N})\in _{H}\mathcal{MHYD}^{H}(\a,\b)$, where $G(N)= N$
  as vector space, with module and comodule structures
  \begin{eqnarray*}
  &h\rightharpoondown n =& \th(h_{1})\b(h_{2})\cdot n, \\
  &\r(n)=:&n^{(0)}\o n^{(1)}= n_{0}\o n_{1}\om.
  \end{eqnarray*}
  Both $F$ and $G$ act as identities on morphisms.

  One checks that $F$ and $G$ are functors and inverse to each other. \hfill $\blacksquare$\\

  {\bf Proposition 4.7.} Let $\a,\b,\g\in {\sl Aut}_{mHH}(H)$.
  The categories $_{H}\mathcal{MHYD}^{H}(\a\b,\g\b)$
and $_{H}\mathcal{MHYD}^{H}(\a,\g)$ are isomorphic.
\\

{\bf Proof.} A pair of inverse functors $(F,G)$ is given as follows.

 If $(M,\xi_{M})\in _{H}\mathcal{MHYD}^{H}(\a\b,\g\b),$
  then $(F(M),\xi_{M})\in _{H}\mathcal{MHYD}^{H}(\a,\g)$, where $F(M)= M$
  as vector space, with structures
  \begin{eqnarray*}
  &h\rightarrow m =& \b^{-1}(h)\cdot m, \\
  &\r(m)=:&m_{<0>}\o m_{<1>}=m_{0}\o m_{1}.
  \end{eqnarray*}

   If $(N,\xi_{N})\in _{H}\mathcal{MHYD}^{H}(\a,\g),$
  then $(G(N),\xi_{N})\in _{H}\mathcal{MHYD}^{H}(\a\b,\g\b)$, where $G(N)= N$
  as vector space, with module and comodule structures
  \begin{eqnarray*}
  &h\rightharpoondown n =& \b(h)\cdot n, \\
  &\r(n)=:&n^{(0)}\o n^{(1)}= n_{0}\o n_{1}.
  \end{eqnarray*}
  Both $F$ and $G$ act on morphisms as identities.

  We can check that $F$ and $G$ are functors and inverse to each other.
  This completes the proof. \hfill $\blacksquare$
  \\

  {\bf Corollary 4.8.} For all $\a,\b\in {\sl Aut}_{mHH}(H)$, we have
  isomorphisms of categories:
  \begin{eqnarray*}
  &&_{H}\mathcal{MHYD}^{H}(\a,\b)\simeq _{H}\mathcal{MHYD}^{H}(\a\b^{-1},id),
  \ \ _{H}\mathcal{MHYD}^{H}(\a,\a)\simeq _{H}\mathcal{MHYD}^{H},\\
  && _{H}\mathcal{MHYD}^{H}(\a,id)\simeq _{H}\mathcal{MHYD}^{H}(id,\a^{-1}),
  \ \ _{H}\mathcal{MHYD}^{H}(id,\b)\simeq _{H}\mathcal{MHYD}^{H}(\a^{-1},id).
  \end{eqnarray*}

  Let again $\a,\b\in  {\sl Aut}_{mHH}(H)$ such that there exist
  $(\th,\om)$ a pair in involution corresponding to $(\a,\b)$, and assume
  that $(H,\xi_{H})$ is finite dimensional. Then we know that $_{H}\mathcal{MHYD}^{H}(\a,\b)
 \simeq  _{H^{*}\bowtie H(\a,\b)}\mathcal{M}$,
 $_{H}\mathcal{MHYD}^{H}\simeq  _{D(H)}\mathcal{M}$
 (\cite{CZ2014}, Proposition 4.3), and the isomorphism $_{H}\mathcal{MHYD}^{H}(\a,\b)\simeq _{H}\mathcal{MHYD}^{H}$ constructed in the theorem is induced by
  a monoidal Hom-algebra isomorphism as follows.
  \\

  {\bf Corollary 4.9.} $(H^{*}\bowtie H(\a,\b),\xi_{H}^{*-1}\o \xi_{H})
  \simeq (D(H), \xi_{H}^{*-1}\o \xi_{H})$ as monoidal Hom-algebras,
   given by
 \begin{eqnarray*}
  D(H)\rightarrow H^{*}\bowtie H(\a,\b),&&
  f\o h\mapsto \om^{-1}\rightharpoonup f\bowtie \th(\b^{-1}(S(h_{1})))\b^{-1}(h_{2}),\\
  H^{*}\bowtie H(\a,\b)\rightarrow D(H),&&
  f\bowtie h\mapsto \om\rightharpoonup f\o \th(h_{1})\b(h_{2}).
  \end{eqnarray*}
  for all $h\in H, f\in H^{*},$ and a group-like element $\om\in H$.
\\

  Finally, we consider some special cases, which are shown in Example 4.3.
  Similar to the cases in Staic \cite{S2007}, we give the following two propositions.

  We define the modular pair $(\om,\th)$ in monoidal Hom-Hopf algebra $(H,\xi_{H})$,
  i.e., $\th$ is an algebra map $H\rightarrow k$ and
  $\om\in(H,\xi_{H})$ is a group-like element
  satisfying $\th(\om)=1.$ Defining an endomorphism $\widetilde{S}$ of $(H,\xi_{H})$
  by $\widetilde{S}(h)=S(h_{1})\th(h_{2})$ for all $h\in H$, then $(\om, \th)$ is called
  a modular pair in involution if $\widetilde{S}^{2}(h)=\om^{-1}(h\om).$
  \\

  {\bf Proposition 4.10.} Let $(H,\xi_{H})$ be a monoidal Hom-Hopf algebra,
  $(\om,\th)$ a modular pair in involution and $(M,\xi_{M})$ a left-right
  anti-Yetter-Drinfeld Hom-module. If we define a new action of $H$ on $M$ as:
  $$h\rightharpoondown m = \th(S(h_{1}))h_{2}\cdot m,$$
  and a new coaction as follows:
  $$\r(m)=m_{<0>}\o m_{<1>}= m_{0}\o m_{1}\om^{-1},$$
 then $(M,\rightharpoondown, \r)$ is a left-right Yetter-Drinfeld Hom-module.
 \\

  {\bf Proof.} First, since $\th: H\rightarrow k$ is an algebra morphism
  and $\om$ is a group-like element, the module and comdule
  structures are given by above formulas.

  We denote the involution inverse of $\th$ by $\th^{-1}$.
  From $\widetilde{S}^{2}(h)=\om^{-1}(h\om)$, we can get
  $\th (S(h_{11}))S^{2}(h_{12})\th(h_{2}) = \om^{-1}(h\om)$
  and $\th^{-1}(h_{1})S(h_{21})\th(h_{22}) = \om^{-1}(S^{-1}(h)\om):$
\begin{eqnarray*}
  &&(h\rightharpoondown m)_{<0>}\o (h\rightharpoondown m)_{<1>}\\
  &=& \th^{-1}(h_{1})(h_{2}\cdot m)_{0}\o (h_{2}\cdot m)_{1}\om^{-1} \\
  &=& \th^{-1}(h_{1})\xi_{H}(h_{221})\cdot m_{0}
  \o ((h_{222} \xi_{H}^{-1}(m_{1}))S(h_{21}))\om^{-1} \\
  &=& h_{21}\cdot m_{0}\o ((\xi_{H}^{-1}(h_{22}) \xi_{H}^{-1}(m_{1}))
  \th^{-1}(\xi_{H}(h_{11}))S(h_{12}))\om^{-1} \\
  &=& h_{21}\cdot m_{0}\o ((\xi_{H}^{-1}(h_{22}) \xi_{H}^{-1}(m_{1}))
  ((\om^{-1}\om)\th^{-1}(\xi_{H}(h_{11}))S(h_{121})\th(h_{1221})\th^{-1}(h_{1222})))\om^{-1} \\
  &=& h_{21}\cdot m_{<0>}\o ((h_{22}\xi_{H}^{-1}(m_{<1>}))
  ((\om^{-1}\th^{-1}(\xi_{H}^{2}(h_{111}))S\xi_{H}(h_{1121})\th(h_{1122}))\om^{-1}))\\
  && \th^{-1}(\xi_{H}^{-2}(h_{12})) \\
  &=& h_{21}\cdot m_{<0>}\o (h_{22}\xi_{H}^{-1}(m_{<1>}))
  S^{-1}\xi_{H}(h_{11})\th^{-1}(\xi_{H}^{-2}(h_{12})) \\
  &=& \th^{-1}(\xi_{H}(h_{211}))\xi_{H}(h_{212})\cdot m_{<0>}
  \o (h_{22}\xi_{H}^{-1}(m_{<1>}))S^{-1}(h_{1}) \\
  &=& \xi_{H}(h_{21})\rightharpoondown m_{<0>}\o (h_{22}\xi_{H}^{-1}(m_{<1>}))S^{-1}(h_{1})
  \end{eqnarray*}
  This means that $(M,\rightharpoondown, \r)$ is a left-right Yetter-Drinfeld Hom-module.
  \hfill $\blacksquare$\\

  By Example 4.4 and Remark 1.2.7 (3), we have the following proposition.
\\

{\bf Proposition 4.11.} For any integer numbers $m$ and $n$, if $(M,\xi_{M})$
is a left-right $m$-Yetter-Drinfeld Hom-module and $(N,\xi_{N})$ is an
$n$-Yetter-Drinfeld Hom-module, then $(M\o N,\xi_{M}\o\xi_{N})$ is
a left-right $m+n$-Yetter-Drinfeld Hom-module with
module structure and comodule structure as follows
\begin{eqnarray*}
  h\c (m \o n) &=& S^{2n} (h_{1})\c m \o h_{2}\c n,\\
 m\o n &\mapsto & (m_{0}\o n_{0})\o n_{1}m_{1}.
  \end{eqnarray*}
  for all $m\in M,n\in N$ and $h\in H.$
  \\

  Let $\mathcal{ZMHYD}(H)$ be the disjoint union of all categories
  $_{H}\mathcal{MHYD}^{H}(S^{2n},id)$ of left-right $n$-Yetter-Drinfeld
  Hom-modules with $n\in \mathbb{Z}$, the set of integer numbers.
  Then by Theorem 3.7 in \cite{YW2014} and Proposition 4.11,
  the following corollary is a generalization of the main result
   in Staic \cite{S2007}.
   \\

   {\bf Corollary 4.12.} $\mathcal{ZMHYD}(H)$ is a braided $T$-category
   over $\mathbb{Z}$.
   \\

{\bf Example 4.13.} Let $A =\<a\>$ be a cyclic group of order $n$,
   and $Aut(A)=\{\s_{t}: \s_{t}(a)=a^{t},0<t<n, (t,n)=1, t\in \mathbb{Z}\}$.
   Then $(k[A],\xi_{k[A]})$ is a monoidal Hom-Hopf algebra with structure
   given by
   \begin{eqnarray*}
   a^{i}\circ a^{j}=\xi_{k[A]}^{-1}(a^{i}a^{j}), &&
   \Delta(a^{i})=\xi_{k[A]}^{-1}(a^{i})\o \xi_{k[A]}^{-1}(a^{i}),\\
   \varepsilon(a^{i})=1_{k}, &&  S(a^{i})= a^{-i},
     \end{eqnarray*}
   for all $i,j\in \mathbb{Z}$.

   First, $Aut_{mHH}(k[A])=Aut(A)$.
   Let $(H,\xi_{H})=(k[A], \xi_{k[A]}=\sigma_{2})$
    is a monoidal Hom-Hopf algebra given by
  \begin{eqnarray*}
  a^{i}\circ a^{j}= \sigma_{2}^{-1}(a^{i}a^{i})= a^{i+j-2}, &&
   \Delta(a^{i})=\sigma_{2}^{-1}(a^{i})\o \sigma_{2}^{-1}(a^{i})=a^{i-2}\o a^{i-2},\\
   \varepsilon(a^{i})=1_{k}, && S(a^{i})= a^{-i}.
  \end{eqnarray*}
   for all $i,j\in \mathbb{Z}$. It is easy to check that
   $$ S^{2n}(a^{i})=a^{i},$$
   for all $n\in \mathbb{Z}.$
   \\

   Let $\mathcal{ZMHYD}(k[A])$ be the disjoint union of all categories
  $_{k[A]}\mathcal{MHYD}^{k[A]}(S^{2n},id)$ of left-right $n-\mathcal{MHYD}$
  with $n\in \mathbb{Z}$.
 \\

 Let $(M,\xi_{M})$ be an $m-\mathcal{MHYD}$-module
  and $(N,\xi_{N})$ be an $n-\mathcal{MHYD}$-module, for all $m,n\in \mathbb{Z}$.
  Then $(M\o N,\xi_{M}\o\xi_{N})$ is $m+n-\mathcal{MHYD}$-module with structures
  as follows:
  \begin{eqnarray*}
   a^{i}\cdot (x\o y) &= &S^{2n}(a^{i-2})\cdot x\o a^{i-2}\cdot y,\\
 (x\o y)&\mapsto &(x_{0}\o y_{0})\o y_{1}x_{1},
    \end{eqnarray*}
    for all $x\in M, y\in N,a^{i}\in k[A], n\in \mathbb{Z}.$

  On $^{(S^{2m},id)}N = N$, there is an action $\unrhd$ given by
  $$ a^{i}\unrhd y = S^{-2m}(a^{i})\cdot y,$$
  and a coaction $\r_{r}$ defined by
  $$y\mapsto y_{0}\o S^{2m}(y_{1}),$$
  $ y\in N,a^{i}\in k[A], m\in \mathbb{Z}.$

   Let $(M,\xi_{M})$ be an $m-\mathcal{MHYD}$-module
  and $(N,\xi_{N})$ be an $n-\mathcal{MHYD}$-module, for all $m,n\in \mathbb{Z}$.
  Then the braiding
  $$c_{M,N}: M\o N \rightarrow ^{(S^{2m},id)}N  \o M$$
  is given by
  $$ c_{M,N}(x\o y)= \xi_{N}(y_{0})\o y_{1}\cdot \xi_{M}^{-1}(x),$$
  for all $x\in M, y\in N, m\in \mathbb{Z}.$

   Then by Corollary 4.12, $\mathcal{ZMHYD}(k[A])$ is
   a new braided $T$-category over $\mathbb{Z}.$

\section*{ACKNOWLEDGEMENTS}

 This work was supported by the NSF of China (No. 11371088) and the NSF of Jiangsu Province (No. BK2012736).

\end{document}